\documentclass[11pt]{article}
\usepackage{epsfig}

\def\be{\begin{equation}}
\def\ee{\end{equation}}
\def\bea{\begin{eqnarray}}
\def\eea{\end{eqnarray}}
\def\bes{\begin{eqnarray*}}
\def\ees{\end{eqnarray*}}

\def\nn{\nonumber}
\def\lb{\label}
\def\bs{\setminus}

\def\T{{\cal T}}
\def\H{{\cal H}}
\def\R{{\bf R}}
\def\C{{\bf C}}
\def\Z{{\bf Z}}

\def\N{{\bf N}}
\def\U{{\bf U}}

\def\Q{{\bf Q}}
\def\T{{\bf T}}

\def\Sg{{\Sigma}}

\def\aa{{\alpha}}
\def\bb{{\beta}}
\def\ga{{\gamma}}

\def\th{{\theta}}

\def\om{{\omega}}
\def\Om{{\Omega}}
\def\ep{{\epsilon}}
\def\lm{{\lambda}}
\def\Lm{{\Lambda}}
\def\dl{{\delta}}

\def\sg{{\sigma}}
\def\dm{{\diamond}}
\def\Sg{{\Sigma}}
\def\vf{{\varphi}}

\def\<{{\langle}}
\def\>{{\rangle}}
\def\T{{\cal T}}

\def\P{{\cal P}}

\def\Nn{{\cal N}}

\def\mul{{\rm mul}}

\def\crit{{\rm crit}}

\def\Sp{{\rm Sp}}

\def\mod{{\rm mod}}

\def\wtd#1{\widetilde{#1}}

\def\hb{\vrule height0.18cm width0.14cm $\,$}

\title{Three closed characteristics on non-degenerate star-shaped hypersurfaces in $\R^{6}$}

\author{Huagui Duan $^{1}$,\thanks{Partially supported by National Key R\&D Program of China (2020YFA0713300)
and NSFC (12271268, 11671215 and 11790271) and the Fundamental Research Funds for the Central Universities. E-mail: duanhg@nankai.edu.cn.}
\quad Hui Liu $^{2}$,\thanks{Partially supported by NSFC (Nos. 12371195, 12022111) and the Fundamental
Research Funds for the Central Universities (No. 2042023kf0207). E-mail: huiliu00031514@whu.edu.cn.}
\quad Yiming Long $^{3}$,\thanks{Partially supported by National Key R$\&$D Program of China (2020YFA0713300), NSFC
(11131004, 11671215 and 11790271), LPMC of Ministry of Education of China, Nankai University, Wenzhong
Foundation, and Nankai Zhide Foundation. E-mail: longym@nankai.edu.cn.}\quad Zihao Qi $^{1}$, \quad Wei Wang $^{4}$\thanks{Partially
supported by NSFC (12025101 and 11431001). E-mail: wangwei@math.pku.edu.cn.}\\
$^{1}$ School of Mathematical Sciences and LPMC, Nankai University, Tianjin 300071\\
$^{2}$ School of Mathematics and Statistics, Wuhan University, Wuhan 430072, Hubei\\
$^{3}$ Chern Institute of Mathematics and LPMC, Nankai University, Tianjin 300071\\
$^{4}$ KLPAM and School of Mathematical Sciences, Peking University, Beijing 100871 \\
The People's Republic of China \\ }

\begin{document}

\maketitle

\begin{abstract}
{\it In this paper, we prove that for every non-degenerate $C^3$ compact star-shaped hypersurface
$\Sigma$ in $\R^{6}$ which carries no prime closed characteristic of Maslov-type index $0$ or no prime closed characteristic of Maslov-type index
$-1$, there exist at least three prime closed characteristics on $\Sigma$.}
\end{abstract}



\renewcommand{\theequation}{\thesection.\arabic{equation}}
\renewcommand{\thefigure}{\thesection.\arabic{figure}}

\setcounter{figure}{0}
\setcounter{equation}{0}
\section{Introduction and main results}

Let $\Sigma$ be a $C^3$ compact hypersurface in $\R^{2n}$ which is strictly star-shaped with respect to the origin, i.e.,
the tangent hyperplane at any $x\in\Sigma$ does not intersect the origin. We denote the set of all such
hypersurfaces by $\H_{st}(2n)$, and denote by $\H_{con}(2n)$ the subset of $\H_{st}(2n)$ which consists of all
strictly convex hypersurfaces. We consider closed characteristics $(\tau, y)$ on $\Sigma$, which are solutions
of the following problem
\be
\left\{\matrix{\dot{y}=JN_\Sigma(y), \cr
               y(\tau)=y(0), \cr }\right. \lb{1.1}\ee
where $J=\left(\matrix{0 &-I_n\cr
        I_n  & 0\cr}\right)$, $I_n$ is the identity matrix in $\R^n$, $\tau>0$, $N_\Sigma(y)$ is the outward
normal vector of $\Sigma$ at $y$ normalized by the condition $N_\Sigma(y)\cdot y=1$. Here $a\cdot b$ denotes
the standard inner product of $a, b\in\R^{2n}$. A closed characteristic $(\tau, y)$ is {\it prime}, if $\tau$
is the minimal period of $y$. Two closed characteristics $(\tau, y)$ and $(\sigma, z)$ are {\it geometrically
distinct}, if $y(\R)\not= z(\R)$. We denote by $\T(\Sigma)$ the set of geometrically distinct
closed characteristics $(\tau, y)$ on $\Sigma$. A closed characteristic
$(\tau, y)$ is {\it non-degenerate} if $1$ is a Floquet multiplier of $y$ of precisely algebraic multiplicity
$2$; {\it hyperbolic} if $1$ is a double Floquet multiplier of it and all the other Floquet multipliers
are not on ${\bf U}=\{z\in {\bf C}\mid |z|=1\}$, i.e., the unit circle in the complex plane; {\it elliptic}
if all the Floquet multipliers of $y$ are on ${\bf U}$; {\it irrationally elliptic} if $1$ is its double Floquet
multipliers, and the other $(2n-2)$ locating on the unit circle with rotation angles being irrational multiples
of $\pi$. We call a $\Sigma\in \H_{st}(2n)$ {\it non-degenerate} if all closed characteristics and their iterates
on $\Sigma$ are non-degenerate.

There is a long standing conjecture on the number of closed characteristics on compact convex hypersurfaces in
$\R^{2n}$:
\be \,^{\#}\T(\Sg)\ge n, \qquad \forall \; \Sg\in\H_{con}(2n). \lb{1.2}\ee

In 1978, P. Rabinowitz in \cite{Rab1} proved
$^\#\T(\Sg)\ge 1$ for any $\Sg\in\H_{st}(2n)$, A. Weinstein in \cite{Wei1}
proved $^\#\T(\Sg)\ge 1$ for any $\Sg\in\H_{con}(2n)$ independently.
When $n\ge 2$,  in 1987-1988, I. Ekeland-L. Lassoued, I. Ekeland-H. Hofer, and A. Szulkin (cf. \cite{EkL1},
\cite{EkH1}, \cite{Szu1}) proved
$$ \,^{\#}\T(\Sg)\ge 2, \qquad \forall\,\Sg\in\H_{con}(2n). $$

In \cite{LoZ} of 2002, Y. Long and C. Zhu further proved
\bea\;^{\#}\T(\Sg)\ge\left[\frac{n}{2}\right]+1, \qquad \forall\, \Sg\in \H_{con}(2n). \nn\eea
In particular, if all the prime closed characteristics on $\Sg$ are non-degenerate, then $\,^{\#}\T(\Sg)\ge n$
(cf. Theorem 1.1 and Corollary 1.1 of \cite{LoZ}). In \cite{WHL} of 2007, W. Wang, X. Hu and Y. Long proved
$\,^{\#}\T(\Sg)\ge 3$ for every $\Sg\in\H_{con}(6)$.  In \cite{Wan1} of 2016, W. Wang proved
$\,^{\#}\T(\Sg)\ge \left[\frac{n+1}{2}\right]+1$ for every $\Sg\in\H_{con}(2n)$. In \cite{Wan2} of 2016, W.
Wang proved $\,^{\#}\T(\Sg)\ge 4$ for every $\Sg\in\H_{con}(8)$.

Note that every contact form supporting the standard contact structure on $M=S^{2n-1}$ arises from embeddings of $M$
into $\R^{2n}$ as a strictly star-shaped hypersurface enclosing the origin, it is conjectured that in fact the
conjecture (\ref{1.2}) holds for any $\Sg\in\H_{st}(2n)$ too. For star-shaped case, \cite{Gir1} of 1984 and
\cite{BLMR} of 1985 show that $\;^{\#}\T(\Sg)\ge n$ for $\Sg\in\H_{st}(2n)$ under some pinching conditions.
In \cite{Vit2} of 1989, C. Viterbo proved a generic existence result for infinitely many closed
characteristics on star-shaped hypersurfaces in ${\bf R}^{4n}$. In \cite{HuL} of 2002, X. Hu and Y. Long proved that
$\;^{\#}\T(\Sg)\ge 2$ for any non-degenerate $\Sg\in \H_{st}(2n)$. In \cite{HWZ2} of 2003, H. Hofer, K. Wysocki, and
E. Zehnder proved that $\,^{\#}\T(\Sg)=2$ or $\infty$ holds for every non-degenerate $\Sg\in\H_{st}(4)$ provided that
the stable and unstable manifolds of every hyperbolic closed orbit on $\Sg$ intersect transversally. This
condition was removed by D. Cristofaro-Gardiner, M. Hutchings and D. Pomerleano \cite{CGHP} of 2019.
Recently, D. Cristofaro-Gardiner, U. Hryniewicz, M. Hutchings and H. Liu in \cite{CGHHL} proved that $\,^{\#}\T(\Sg)=2$ or $\infty$ holds for every $\Sg\in\H_{st}(4)$. In \cite{CGH1} of 2016, D. Cristofaro-Gardiner and M. Hutchings proved that $\;^{\#}\T(\Sg)\ge 2$ for any contact three
manifold $\Sg$. Various proofs of this result for a star-shaped hypersurface can be found in \cite{GHHM}, \cite{LLo}
and \cite{GiG1}. In \cite{GuK} of 2016, J. Gutt and J. Kang proved $\,^{\#}\T(\Sg)\ge n$ for every  non-degenerate
$\Sg\in \H_{st}(2n)$ if every closed characteristic on $\Sg$ possesses Conley-Zehnder index at least $(n-1)$. In
\cite{DLLW1} of 2018, H. Duan, H. Liu, Y. Long and W. Wang proved $\,^{\#}\T(\Sg)\ge n$ for every index perfect
non-degenerate $\Sg\in \H_{st}(2n)$, and there exist at least $n$ (or $(n-1)$) non-hyperbolic closed characteristics
when $n$ is even (or odd). Here $\Sg\in \H_{st}(2n)$ is {\it index perfect} if it carries only finitely many
geometrically distinct prime closed characteristics, and every prime closed characteristic $(\tau,y)$ on $\Sigma$
possesses positive mean index and whose Maslov-type index $i(y, m)$ of its $m$-th iterate satisfies $i(y, m)\not= -1$
when $n$ is even, and $i(y, m)\not\in \{-2,-1,0\}$ when $n$ is odd for all $m\in\N$. Then later V. Ginzburg, B. G\"{u}rel
and L. Macarini in \cite{GGM} extended these results to prequantization bundles which covered many known results for
closed geodesics on Finsler manifolds.

This paper is mainly devoted to studying the above conjecture on non-degenerate $\Sigma$ in $\R^{6}$. In order to obtain
the existence of three prime closed characteristics on such $\Sigma$, we remove the assumption on the positive mean
index, and relax the index condition $i(y, m)\not\in \{-2,-1,0\}$ for all $m\in\N$ in \cite{DLLW1} to $i(y,1)\neq 0$
or $i(y,1)\neq -1$ for every prime closed characteristic $y$ as follows.

\medskip

{\bf Theorem 1.1.} {\it For every non-degenerate $C^3$ compact star-shaped hypersurface $\Sigma$ in $\R^{6}$, there
exist at least three prime closed characteristics, provided that one of the following conditions holds: (i) there
exists no prime closed characteristic on $\Sigma$ possessing Maslov-type index $0$; (ii) there exists no prime closed
characteristic on $\Sigma$ possessing Maslov-type index $-1$.}

\medskip

Our main idea is the following. For every non-degenerate $\Sg\in \H_{st}(6)$, there holds $\;^{\#}\T(\Sg)\ge 2$
by \cite{HuL}. So firstly we assume that $\;^{\#}\T(\Sg)\ge 2$. Then through some precise analysis of these two
prime closed characteristics $y_1$ and $y_2$, we obtain many properties about their indices and classifications.
For example, (i) both of their mean indices are positive (see Lemma 4.2 below), where we follow some ideas from
\cite{Vit1} and \cite{Vit2}; (ii) These two closed characteristics must belong to one of two different classes
(see (\ref{4.49}) and (\ref{4.50}) in Lemma 4.6 below). With the help of these properties, we conclude that
their Maslov-type indices must satisfy $\{i(y_1,1),i(y_2,1)\}=\{0,-1\}$ (see (\ref{4.56})-(\ref{4.57}) in Lemma
4.9 and (\ref{4.67})-(\ref{4.68}) in Lemma 4.10 below), where we especially use a variant of the generalized
common index jump theorem from \cite{DuQ} (see Theorems 3.2 and 3.3 below) and some corresponding arguments.
Note that in Theorem 1.1 we assume that $\Sigma$ does not possess any prime closed characteristic with either
Maslov-type index of zero or $-1$, thus we obtain a contradiction through the above arguments, which implies
that $\;^{\#}\T(\Sg)\ge 3$ as claimed in Theorem 1.1.

\medskip

In this paper, let $\N$, $\N_0$, $\Z$, $\Q$, $\R$,  $\C$ and $\R^+$ denote
the sets of natural integers, non-negative integers, integers,
rational numbers, real numbers, complex numbers and positive real numbers
respectively. Denote by $a\cdot b$ and $|a|$ the standard inner
product and norm in $\R^{2n}$. Denote by $\langle\cdot,\cdot\rangle$
and $\|\cdot\|$ the standard $L^2$-inner product and $L^2$-norm. For
an $S^1$-space $X$, we denote by $X_{S^1}$ the homotopy quotient of
$X$ module the $S^1$-action, i.e., $X_{S^1}=S^\infty\times_{S^1}X$.
We define the functions \bea \left\{\matrix{[a]=\max\{k\in\Z\,|\,k\le
a\}, & E(a)=\min\{k\in\Z\,|\,k\ge a\} , \cr
                   \varphi(a)=E(a)-[a],  & \{a\}=a-[a]. \cr}\right. \lb{1.3}\eea
Specially, $\varphi(a)=0$ if $ a\in\Z\,$, and $\varphi(a)=1$ if
$a\notin\Z\,$. In this paper we use only $\Q$-coefficients for all
homological modules. For a $\Z_m$-space pair $(A, B)$, let
$H_{\ast}(A, B)^{\pm\Z_m}= \{\sigma\in H_{\ast}(A,
B)\,|\,L_{\ast}\sigma=\pm \sigma\}$, where $L$ is a generator of the
$\Z_m$-action.

\setcounter{figure}{0}
\setcounter{equation}{0}
\section{Critical point theory for closed characteristics on compact star-shaped hypersurfaces in $\R^{2n}$}

In this section, we review briefly  the critical point theory for
closed characteristics on $\Sg\in\H_{st}(2n)$ developed in \cite{Vit1} and \cite{LLW} which will be needed in Section 4. All
the details of proofs can be found in \cite{Vit1} and \cite{LLW}.
Now we fix a $\Sg\in\H_{st}(2n)$ and assume the following condition:

\medskip

(F) {\bf There exist only finitely many geometrically distinct prime closed characteristics\\
$\qquad\qquad \{(\tau_j, y_j)\}_{1\le j\le k}$ on $\Sigma$. }

\medskip

Let $\hat{\tau}=\inf_{1\leq j\leq k}{\tau_j}$ and $T$ be a fixed positive constant. Then by Section 2 of
\cite{LLW}, for any $a>\frac{\hat{\tau}}{T}$, we can construct a function $\varphi_a\in C^{\infty}({\bf R}, {\bf R}^+)$
which has $0$ as its unique critical point in $[0, +\infty)$. Moreover, $\frac{\varphi^{\prime}(t)}{t}$ is strictly
decreasing for $t>0$ together with $\varphi(0)=0=\varphi^{\prime}(0)$ and
$\varphi^{\prime\prime}(0)=1=\lim_{t\rightarrow 0^+}\frac{\varphi^{\prime}(t)}{t}$. More precisely, we
define $\varphi_a$ and the Hamiltonian function $\wtd{H}_a(x)=a\vf_a(j(x))$ via Lemma 2.2 and Lemma 2.4 in \cite{LLW}.
Then we can define a new Hamiltonian function $H_a$ via Proposition 2.5 of \cite{LLW} and consider the fixed period
problem
\be \left\{\matrix{\dot{x}(t) = JH_a^\prime(x(t)), \cr
                         x(0) = x(T).  \cr }\right. \lb{2.1}\ee
Then $H_a\in C^{3}({\bf R}^{2n} \setminus\{0\},{\bf R})\cap C^{1}({\bf R}^{2n},{\bf R})$.
Solutions of (\ref{2.1}) are $x\equiv 0$ and $x=\rho y(\tau t/T)$ with
$\frac{\vf_a^\prime(\rho)}{\rho}=\frac{\tau}{aT}$, where $(\tau, y)$ is a solution of (\ref{1.1}). In particular,
non-zero solutions of (\ref{2.1}) are in one to one correspondence with solutions of (\ref{1.1}) with period
$\tau<aT$.

For any $a>\frac{\hat{\tau}}{T}$, choose some large constant $K=K(a)$ such that
$$ H_{a,K}(x) = H_a(x)+\frac{1}{2}K|x|^2   $$
is a strictly convex function, i.e.,
$$ (\nabla H_{a, K}(x)-\nabla H_{a, K}(y), x-y) \geq \frac{\ep}{2}|x-y|^2,  $$
for all $x, y\in {\bf R}^{2n}$ and some $\ep>0$. Let $H_{a,K}^*$ be the Fenchel dual of $H_{a,K}$
defined by
$$  H_{a,K}^\ast (y) = \sup\{x\cdot y-H_{a,K}(x)\;|\; x\in \R^{2n}\}.   $$
The dual action functional on $X=W^{1, 2}({\bf R}/{T {\bf Z}}, {\bf R}^{2n})$ is defined by
\be F_{a,K}(x) = \int_0^T{\left[\frac{1}{2}(J\dot{x}-K x,x)+H_{a,K}^*(-J\dot{x}+K x)\right]dt}. \lb{2.2}\ee
It can be proved that  $F_{a,K}\in C^{1,1}(X, \R)$ and for $KT\not\in 2\pi{\bf Z}$, $F_{a,K}$ satisfies the
Palais-Smale condition,  $x$ is a critical point of $F_{a, K}$ if and only if it is a solution of (\ref{2.1}). Moreover,
for every critical point $x_a\neq 0$ of $F_{a, K}$, we have $F_{a, K}(x_a)<0$ and it is independent of $K$.

For $KT\notin 2\pi{\bf Z}$, the map $x\mapsto -J\dot{x}+Kx$ is an  isomorphism between the Hilbert spaces
$X=W^{1, 2}({\bf R}/({T {\bf Z}}); {\bf R}^{2n})$ and $E=L^{2}({\bf R}/(T {\bf Z}),{\bf R}^{2n})$. Denote its inverse
by $M_K$ and define the functional
\be \Psi_{a,K}(u)=\int_0^T{\left[-\frac{1}{2}(M_{K}u, u)+H_{a,K}^*(u)\right]dt}, \qquad \forall\,u\in E. \lb{2.3}\ee
Then $x\in X$ is a critical point of $F_{a,K}$ if and only if $u=-J\dot{x}+Kx$ is a critical point of $\Psi_{a, K}$.

Suppose $u$ is a nonzero critical point of $\Psi_{a, K}$,  the formal Hessian of $\Psi_{a, K}$ at $u$ is defined by
$$ Q_{a,K}(v)=\int_0^T(-M_K v\cdot v+H_{a,K}^{*\prime\prime}(u)v\cdot v)dt,  $$
which defines an orthogonal splitting $E=E_-\oplus E_0\oplus E_+$ of $E$ into negative, zero and positive subspaces.
The index and nullity of $u$ are defined by $i_K(u)=\dim E_-$ and $\nu_K(u)=\dim E_0$ respectively.
Similarly, we define the index and nullity of $x=M_Ku$ for $F_{a, K}$ and denote them by $i_K(x)$ and
$\nu_K(x)$. Then it follows from  (\ref{2.2}) and (\ref{2.3}) that
\be  i_K(u)=i_K(x),\quad \nu_K(u)=\nu_K(x).  \lb{2.4}\ee
By Lemma 6.4 of \cite{Vit2}, we have
\be  i_K(x) = 2n([KT/{2\pi}]+1)+i(x) \equiv d(K)+i(x),   \lb{2.5}\ee
where the index $i(x)$ does not depend on K, but only on $H_a$. This index was first introduced by C. Viterbo in
\cite{Vit2}.

By the proof of Proposition 2 of \cite{Vit1}, $v\in E$ belongs to the null space of $Q_{a, K}$
if and only if $z=M_K v$ is a solution of the linearized system
$$  \dot{z}(t) = JH_a''(x(t))z(t).  $$
Thus the nullity in (\ref{2.4}) is independent of $K$, and we denote it by $\nu(x)\equiv \nu_K(u)= \nu_K(x)$.

Note that by Proposition 2.11 of \cite{LLW}, the index $i(x)$ and nullity $\nu(x)$ coincide with those defined for
the Hamiltonian $H(x)=j(x)^\alpha$ for all $x\in\R^{2n}$ and some $\aa\in (1,2)$. Especially
$1\le \nu(x)\le 2n-1$ always holds.

Suppose $(\tau, y)$ is a closed characteristic on $\Sigma$, let $a>\tau/T$ and choose $\vf_a$ as defined above
(\ref{2.1}). Determine $\rho$ uniquely by $\frac{\vf_a'(\rho)}{\rho}=\frac{\tau}{aT}$. Let
$x=\rho y(\frac{\tau t}{T})$. Then as in \cite{LLW} we define the index $i(\tau,y)$ and nullity $\nu(\tau,y)$ of
$(\tau,y)$ by
\be i(\tau,y)=i(x), \qquad \nu(\tau,y)=\nu(x),  \lb{2.6}\ee
and the mean index of $(\tau,y)$ by
\be \hat{i}(\tau,y) = \lim_{m\rightarrow\infty}\frac{i(m\tau,y)}{m}.  \lb{2.7}\ee
By Proposition 2.11 of \cite{LLW}, the index and nullity are well defined and independent of the choice of $a$. For
a closed characteristic $(\tau,y)$ on $\Sigma$, we simply denote by $y^m\equiv(m\tau,y)$ the $m$-th iteration of $y$
and by $i(y^m)=i(m\tau,y^m)$ and $\nu(y^m)=\nu(m\tau,y^m)$ the index and nullity of $y^m$ for $m\in\N$ in the proofs
below.

There is a natural $S^1$-action on $X$ or $E$ defined by
$$  \theta\cdot u(t)=u(\theta+t),\quad\forall\, \theta\in S^1, \, t\in\R.  $$
Clearly both of $F_{a, K}$ and $\Psi_{a, K}$ are $S^1$-invariant. For any $\kappa\in\R$, let
\bea
\Lambda_{a, K}^\kappa &=& \{u\in L^{2}({\bf R}/({T {\bf Z}}); {\bf R}^{2n})\;|\;\Psi_{a,K}(u)\le\kappa\},  \nn\\
X_{a, K}^\kappa &=& \{x\in W^{1, 2}({\bf R}/(T {\bf Z}),{\bf R}^{2n})\;|\;F_{a, K}(x)\le\kappa\}.  \nn\eea
For a critical point $u$ of $\Psi_{a, K}$ and  $x=M_K u$ the corresponding critical point of $F_{a, K}$, let
\bea
\Lm_{a,K}(u) &=& \Lm_{a,K}^{\Psi_{a, K}(u)}
   = \{w\in L^{2}(\R/(T\Z), \R^{2n}) \;|\; \Psi_{a, K}(w)\le\Psi_{a,K}(u)\},  \nn\\
X_{a,K}(x) &=& X_{a,K}^{F_{a,K}(x)} = \{y\in W^{1, 2}(\R/(T\Z), \R^{2n}) \;|\; F_{a,K}(y)\le F_{a,K}(x)\}. \nn\eea
Then both sets are $S^1$-invariant. Denote by $\crit(\Psi_{a, K})$ the set of critical points of $\Psi_{a, K}$.
Since $\Psi_{a,K}$ is $S^1$-invariant, if $u\in \crit(\Psi_{a, K})$, then $S^1\cdot u\subset \crit(\Psi_{a, K})$ is
a critical orbit of $\Psi_{a, K}$. Note that by the condition (F), the number of critical orbits of $\Psi_{a, K}$
is finite, similarly for $F_{a, K}$. Hence as usual we can make the following definition.

\medskip

{\bf Definition 2.1.} {\it Suppose $u$ is a nonzero critical point of $\Psi_{a, K}$, and $\Nn$ is an $S^1$-invariant
open neighborhood of $S^1\cdot u$ such that $\crit(\Psi_{a,K})\cap (\Lm_{a,K}(u)\cap \Nn) = S^1\cdot u$.
Then the $S^1$-critical module of $S^1\cdot u$ is defined by
$$ C_{S^1,\; q}(\Psi_{a, K}, \;S^1\cdot u)
    = H_{q}((\Lambda_{a, K}(u)\cap\Nn)_{S^1},\; ((\Lambda_{a,K}(u)\setminus S^1\cdot u)\cap\Nn)_{S^1}). $$
Similarly, we can define $C_{S^1,\; q}(F_{a, K}, \;S^1\cdot x)$.}

\medskip

Fix $a$ and let $u_K\neq 0$ be a critical point of $\Psi_{a, K}$ with multiplicity $\mul(u_K)=1$,
i.e., $u_K$ corresponds to a prime closed characteristic $(\tau, y)\subset\Sigma$.
Then we have $u_K=-J\dot x+Kx$ with $x$
being a solution of (\ref{2.1}) and $x=\rho y(\frac{\tau t}{T})$ with
$\frac{\vf_a^\prime(\rho)}{\rho}=\frac{\tau}{aT}$.
For any $p\in\N$ satisfying $p\sigma<aT$, we choose $K$
such that $pK\notin \frac{2\pi}{T}\Z$, then the $p$-th iteration $u_{pK}^p$ of $u_K$ is given by $-J\dot x^p+pKx^p$,
where $x^p$ is the unique solution of (\ref{2.1}) corresponding to $(p\tau, y)$
and is a critical point of $F_{a, pK}$, i.e.,
$u_{pK}^p$ is the critical point of $\Psi_{a, pK}$ corresponding to $x^p$.

\medskip

{\bf Lemma 2.2.} (cf. Proposition 4.2 and Remark 4.4 of \cite{LLW} ) {\it
Suppose $u_K\neq 0$ be a critical point of $\Psi_{a, K}$ satisfying $\mul(u_K)=1$.
If $u_{pK}^p$ is non-degenerate,
i.e., $\nu_{pK}(u_{pK}^p)=1$, then we have}
\bea C_{S^1,q-d(pK)+d(K)}(F_{a,K},S^1\cdot x^p)
&=& C_{S^1,q}(F_{a,pK},S^1\cdot x^p)=C_{S^1,q}(\Psi_{a,pK},S^1\cdot u^p_{pK}) \nn\\
&=& \left\{\matrix{
     \Q, &\quad {\it if}\;\; q=i_{pK}(u_{pK}^p),\;\;{\it and}\;\;\bb(x^p)=1, \cr
     0, &\quad {\it otherwise}, \cr}\right.  \lb{2.8}\eea
where $\bb(x^p)=(-1)^{i_{pK}(u_{pK}^p)-i_{K}(u_{K})}=(-1)^{i(x^p)-i(x)}$ by (\ref{2.4})-(\ref{2.5}).

\medskip

We have the following mean index identity for closed characteristics on $\Sg$.

\medskip

{\bf Theorem 2.3.} (cf.  Theorem 1.2 of \cite{Vit2} and Theorem 1.1 of \cite{LLW}) {\it Suppose that
$\Sg\in\H_{st}(2n)$ satisfying $^\#\T(\Sg)<+\infty$. Denote by $\{(\tau_j,y_j)\}_{1\le j\le k}$ all the
geometrically distinct prime closed characteristics. Then the following identities hold
\be \sum_{1\le j\le k \atop \hat{i}(y_j)>0}\frac{\hat{\chi}(y_j)}{\hat{i}(y_j)}=\frac{1}{2},\qquad
         \sum_{1\le j\le k \atop \hat{i}(y_j)<0}\frac{\hat{\chi}(y_j)}{\hat{i}(y_j)}=0,  \lb{2.9}\ee
where $\hat{\chi}(y)\in\Q$ is the average Euler characteristic given by Definition 4.8 and Remark 4.9 of \cite{LLW}.

In particular, if all $y^m$'s are non-degenerate for $m\ge 1$, then
\bea \hat{\chi}(y)=\left\{\matrix{
     (-1)^{i(y)}, &\quad {\it if}\;\; i(y^2)-i(y)\in 2\Z, \cr
     \frac{(-1)^{i(y)}}{2}, &\quad {\it otherwise}. \cr}\right.  \lb{2.10}\eea}

\medskip

Let $F_{a, K}$ be a functional defined by (\ref{2.2}) for some $a, K\in\R$ large enough and let $\ep>0$ be
small enough such that $[-\ep, 0)$ contains no critical values of $F_{a, K}$. For $b$ large enough,
The normalized Morse series of $F_{a, K}$ in $ X^{-\ep}\setminus X^{-b}$
is defined by
\be  M_a(t)=\sum_{q\ge 0,\;1\le j\le p} \dim C_{S^1,\;q}(F_{a, K}, \;S^1\cdot v_j)t^{q-d(K)},  \lb{2.11}\ee
where we denote by $\{S^1\cdot v_1, \ldots, S^1\cdot v_p\}$ the critical orbits of $F_{a, K}$ with critical
values less than $-\ep$. The Poincar\'e series of $H_{S^1, *}( X, X^{-\ep})$ is $t^{d(K)}Q_a(t)$, according
to Theorem 5.1 of \cite{LLW}, if we set $Q_a(t)=\sum_{k\in \Z}{q_kt^k}$, then
$$   q_k=0,\qquad\forall\;k\in \mathring {I},  $$
where $I$ is an interval of $\Z$ such that $I \cap [i(\tau, y), i(\tau, y)+\nu(\tau, y)-1]=\emptyset$ for all
closed characteristics $(\tau,\, y)$ on $\Sigma$ with $\tau\ge aT$. Then by Section 6 of \cite{LLW}, we have
$$  M_a(t)-\frac{1}{1-t^2}+Q_a(t) = (1+t)U_a(t),   $$
where $U_a(t)=\sum_{i\in \Z}{u_it^i}$ is a Laurent series with nonnegative coefficients.
If there is no closed characteristic with $\hat{i}=0$, then
\be   M(t)-\frac{1}{1-t^2}=(1+t)U(t),    \lb{2.12}\ee
where $M(t)=\sum_{p\in \Z}{M_pt^p}$ denotes $M_a(t)$ as $a$ tends to infinity. In addition, we also denote
by $b_p$ the coefficient of $t^p$ of $\frac{1}{1-t^2}=\sum_{p\in \Z}{b_pt^p}$, i.e. there holds $b_p=1$,
$\forall\ p\in2\N_0$ and $b_p=0$, $\forall\ p\not\in2\N_0$.

\setcounter{figure}{0}
\setcounter{equation}{0}
\section{Index iteration theory for symplectic paths}

In \cite{Lon1} of 1999, Y. Long established the basic normal form decomposition of symplectic matrices.
Based on this result he further established the precise iteration formulae of the Maslov-type indices for
symplectic paths in \cite{Lon2} of 2000.

As in \cite{Lon2}, denote by
\bea
N_1(\lm, b) &=& \left(\matrix{\lm & b\cr
                                0 & \lm\cr}\right), \qquad {\rm for\;}\lm=\pm 1, \; b\in\R, \nn\\
D(\lm) &=& \left(\matrix{\lm & 0\cr
                      0 & \lm^{-1}\cr}\right), \qquad {\rm for\;}\lm\in\R\bs\{0, \pm 1\}, \nn\\
R(\th) &=& \left(\matrix{\cos\th & -\sin\th \cr
                           \sin\th & \cos\th\cr}\right), \qquad {\rm for\;}\th\in (0,\pi)\cup (\pi,2\pi), \nn\\
N_2(e^{\th\sqrt{-1}}, B) &=& \left(\matrix{ R(\th) & B \cr
                  0 & R(\th)\cr}\right), \qquad {\rm for\;}\th\in (0,\pi)\cup (\pi,2\pi)\;\; {\rm and}\; \nn\\
        && \qquad B=\left(\matrix{b_1 & b_2\cr
                                  b_3 & b_4\cr}\right)\; {\rm with}\; b_j\in\R, \;\;
                                         {\rm and}\;\; b_2\not= b_3. \nn\eea
Here $N_2(e^{\th\sqrt{-1}}, B)$ is non-trivial if $(b_2-b_3)\sin\theta<0$, and trivial
if $(b_2-b_3)\sin\theta>0$.

As in \cite{Lon2} and \cite{Lon3}, the $\diamond$-sum of any two real matrices is defined by
$$ \left(\matrix{A_1 & B_1\cr C_1 & D_1\cr}\right)_{2i\times 2i}\diamond
      \left(\matrix{A_2 & B_2\cr C_2 & D_2\cr}\right)_{2j\times 2j}
=\left(\matrix{A_1 & 0 & B_1 & 0 \cr
                                   0 & A_2 & 0& B_2\cr
                                   C_1 & 0 & D_1 & 0 \cr
                                   0 & C_2 & 0 & D_2}\right). $$
For $M\in\Sp(2n)$, we denote its $k$-copy $\diamond$-sum $M\diamond \cdots \diamond M$ by $M^{\diamond k}$.

For every $M\in\Sp(2n)$, the homotopy set $\Omega(M)$ of $M$ in $\Sp(2n)$ is defined in \cite{Lon2}
(cf. Definition 1.8.5 of \cite{Lon3}) by
$$ \Om(M)=\{N\in\Sp(2n)\,|\,\sg(N)\cap\U=\sg(M)\cap\U\equiv\Gamma\;\mbox{and}
                    \;\nu_{\om}(N)=\nu_{\om}(M),\, \forall\om\in\Gamma\}, $$
where $\sg(M)$ denotes the spectrum of $M$, $\nu_{\om}(M)\equiv\dim_{\C}\ker_{\C}(M-\om I)$ for $\om\in\U$.
The set $\Om^0(M)$  is defined to be the path connected component of $\Om(M)$ containing $M$.

For every $\ga\in\mathcal{P}_\tau(2n)\equiv\{\ga\in C([0,\tau],Sp(2n))\ |\ \ga(0)=I_{2n}\}$, we extend
$\ga(t)$ to $t\in [0,m\tau]$ for every $m\in\N$ by
$$ \ga^m(t)=\ga(t-j\tau)\ga(\tau)^j \qquad \forall\;j\tau\le t\le (j+1)\tau \;\;
               {\rm and}\;\;j=0, 1, \ldots, m-1, $$
as in p.114 of \cite{Lon1}. As in \cite{LoZ} and \cite{Lon3}, we denote the Maslov-type indices of
$\ga^m$ by $(i(\ga,m),\nu(\ga,m))$. Note that these differences of these index and nullity from those defined
in (\ref{2.6}) are given in Theorem 3.2 below.

The following is the precise index iteration formulae for symplectic paths, which is due to Y. Long (cf. Theorems
1.2 and 1.3 of \cite{Lon2} and Theorems 8.2.1 and 8.3.1 in Chapter 8 of \cite{Lon3}).

\medskip

{\bf Theorem 3.1.} {\it Let $\ga\in\P_{\tau}(2n)$. Then there exists a path $f\in C([0,1],\Omega^0(\gamma(\tau))$
such that $f(0)=\gamma(\tau)$ and
\bea f(1)
&=& N_1(1,1)^{\diamond p_-} \diamond I_{2p_0}\diamond N_1(1,-1)^{\diamond p_+}
       \diamond N_1(-1,1)^{\diamond q_-} \diamond (-I_{2q_0})\diamond N_1(-1,-1)^{\diamond q_+}\nn\\
&& \diamond R(\theta_1)\diamond\cdots\diamond R(\theta_r)\diamond N_2(\omega_1, u_1)\diamond\cdots
       \diamond N_2(\omega_{r_*}, u_{r_*}) \nn\\
&& \diamond N_2(\lm_1, v_1)\diamond\cdots\diamond N_2(\lm_{r_0}, v_{r_0})\diamond M_s \nn\eea
where $N_2(\omega_j, u_j)$s are non-trivial and $N_2(\lm_j, v_j)$s are trivial basic normal forms; $\sigma(M_s)\cap U=\emptyset$
with $M_s=D(2)^{\diamond s}$ with an integer $s\ge 0$ or $M_s=D(-2)\diamond D(2)^{\diamond (s-1)}$ with an integer $s\ge 1$;
$p_-$, $p_0$, $p_+$, $q_-$, $q_0$, $q_+$, $r$, $r_*$ and $r_0$ are non-negative integers; $\omega_j=e^{\sqrt{-1}\alpha_j}$,
$\lambda_j=e^{\sqrt{-1}\beta_j}$; $\theta_j$, $\alpha_j$, $\beta_j\in (0, \pi)\cup (\pi, 2\pi)$; these integers and real
numbers are uniquely determined by $\gamma(\tau)$. Then using the functions defined in (\ref{1.3}), we have
\bea i(\gamma, m)
&=& m(i(\gamma,1)+p_-+p_0-r)+2\sum_{j=1}^r E\left(\frac{m\theta_j}{2\pi}\right)-r-p_--p_0\nn\\
& & -\frac{1+(-1)^m}{2}(q_0+q_+)+2\left(\sum_{j=1}^{r_*}\varphi\left(\frac{m\alpha_j}{2\pi}\right)-r_*\right), \nn\eea
\bea \nu(\gamma, m)
&=& \nu(\gamma,1)+\frac{1+(-1)^m}{2}(q_-+2q_0+q_+)+2(r+r_*+r_0)\nn\\
& & -2\left(\sum_{j=1}^{r}\varphi\left(\frac{m\theta_j}{2\pi}\right)+\sum_{j=1}^{r_*}\varphi\left(\frac{m\alpha_j}{2\pi}\right)
    +\sum_{j=1}^{r_0}\varphi\left(\frac{m\beta_j}{2\pi}\right)\right),\nn\\
\hat{i}(\gamma, 1) &=& \lim_{m\to +\infty}\frac{i(\gamma,m)}{m}=i(\gamma, 1)+p_-+p_0-r+\sum_{j=1}^r\frac{\theta_j}{\pi},\nn\eea
where $N_1(1, \pm 1)=
\left(\matrix{ 1 &\pm 1\cr 0 & 1\cr}\right)$, $N_1(-1, \pm 1)=
\left(\matrix{ -1 &\pm 1\cr 0 & -1\cr}\right)$,
$R(\theta)=\left(\matrix{\cos\th &-\sin\th \cr
                         \sin\th & \cos\th \cr}\right)$,
$N_2(\omega, b)=\left(\matrix{R(\th) & b\cr
                              0 & R(\th)\cr}\right)$ with some $\th\in (0,\pi)\cup (\pi,2\pi)$ and
$b=\left(\matrix{b_1 & b_2\cr
                 b_3 & b_4\cr}\right)\in\R^{2\times2}$,
such that $(b_2-b_3)\sin\theta>0$, if $ N_2(\omega, b)$ is trivial; $(b_2-b_3)\sin\theta<0$, if
$N_2(\omega, b)$ is non-trivial.

We have $i(\gamma, 1)$ is odd if $f(1)=N_1(1, 1)$, $I_2$, $N_1(-1, 1)$, $-I_2$,
$N_1(-1, -1)$ and $R(\theta)$; $i(\gamma, 1)$ is even if $f(1)=N_1(1, -1)$ and $N_2(\omega, b)$; $i(\gamma,1)$
can be any integer if $\sigma(f(1))\cap\U = \emptyset$.}

\medskip

The common index jump theorem (cf. Theorem 4.3 of \cite{LoZ}) for
symplectic paths has become one of the main tools in studying the
multiplicity and stability of periodic orbits in Hamiltonian and
symplectic dynamics. Then this theorem has been generalized by
\cite{DLW} and \cite{DLLW2} to the following theorem.

\medskip

{\bf Theorem 3.2.} (cf. Theorem 3.6 of \cite{DLLW2}, {\bf Generalized common index jump theorem})
{\it Let $\gamma_j\in\mathcal{P}_{\tau_j}(2n)$ for $j=1,\cdots,q$ be a finite collection of symplectic paths with nonzero mean
indices $\hat{i}(\ga_j,1)$. Let $M_j=\ga_j(\tau_j)$. Then for any fixed $\bar{m}\in \N$, there exist infinitely many $(q+1)$-tuples
$(N, m_1,\cdots,m_q) \in \N^{q+1}$ such that the following hold for all $1\le j\le q$ and $1\le m\le \bar{m}$,
\bea
\nu(\ga_j,2m_j-m) &=& \nu(\ga_j,2m_j+m) = \nu(\ga_j, m),   \lb{3.1}\\
i(\ga_j,2m_j+m) &=& 2\varrho_j N+i(\ga_j,m),                         \lb{3.2}\\
i(\ga_j,2m_j-m) &=&  2\varrho_j N-i(\ga_j,m)-2(S^+_{M_j}(1)+Q_j(m)),  \lb{3.3}\\
i(\ga_j, 2m_j)&=& 2\varrho_j N -(S^+_{M_j}(1)+C(M_j)-2\Delta_j),     \lb{3.4}\eea
where \bea &&\varrho_j=\left\{\begin{array}{cc}1, &{\rm if}\ \hat{i}(\ga_j,1)>0, \cr
                                     -1, &{\rm if}\  \hat{i}(\ga_j,1)<0, \end{array}\right.\qquad
\Delta_j = \sum_{0<\{m_j\th/\pi\}<\delta}S^-_{M_j}(e^{\sqrt{-1}\th}),\nn\\
&&\ Q_j(m) = \sum_{e^{\sqrt{-1}\th}\in\sg(M_j),\atop \{\frac{m_j\th}{\pi}\}
                   = \{\frac{m\th}{2\pi}\}=0}S^-_{M_j}(e^{\sqrt{-1}\th}). \lb{3.5}\eea
Moreover we have
\bea \min\left\{\{\frac{m_j\theta}{\pi}\},1-\{\frac{m_j\theta}{\pi}\}\right\}<\dl,\lb{3.6}\eea
whenever $e^{\sqrt{-1}\theta}\in\sigma(M_j)$ and $\dl$ can be chosen as small as we want. More precisely,
by (3.17) in \cite{DLLW2} and (4.40), (4.41) in \cite{LoZ} , we have
\bea m_j=\left(\left[\frac{N}{M|\hat i(\gamma_j, 1)|}\right]+\chi_j\right)M,\quad\forall\  1\le j\le q,\lb{3.7}\eea
where $\chi_j=0$ or $1$ for $1\le j\le q$ and $\frac{M\theta}{\pi}\in\Z$
whenever $e^{\sqrt{-1}\theta}\in\sigma(M_j)$ and $\frac{\theta}{\pi}\in\Q$
for some $1\le j\le q$.  Furthermore, given $M_0$, from the proof of Theorem 4.1 of \cite{LoZ}, we may
further require $N$ to be the multiple of $M_0$, i.e., $M_0|N$.}

\medskip

Most recently, Duan and Qi have obtained a variant of the generalized common index jump theorem in \cite{DuQ} as follows.

\medskip

{\bf Theorem 3.3.} (cf. Theorem 3.8 of \cite{DuQ})  {\it Let $\gamma_j\in\mathcal{P}_{\tau_j}(2n)$ for
$j=1,\cdots,q$ be a finite collection of symplectic paths with nonzero mean indices $\hat{i}(\ga_j,1)$.
Let $M_j=\ga_j(\tau_j)$. Fix a positive integer $\hat{p}$ and then let $N$ be an integer satisfying
(\ref{3.1})-(\ref{3.7}).

Then for $\hat{N}=\hat{p}N$, there exists the corresponding $q$-tuples $(\hat{m}_1,\cdots,\hat{m}_q) \in \N^q$
such that the following hold for all $1\le j\le q$ and $1\le m\le \bar{m}$
\bea
\nu(\ga_j,2\hat{m}_j-m) &=& \nu(\ga_j,2\hat{m}_j+m) = \nu(\ga_j, m),   \lb{3.8}\\
i(\ga_j,2\hat{m}_j+m) &=& 2\varrho_j \hat{p} N+i(\ga_j,m),                         \lb{3.9}\\
i(\ga_j,2\hat{m}_j-m) &=&  2\varrho_j \hat{p} N-i(\ga_j,m)-2(S^+_{M_j}(1)+Q_j(m)),  \lb{3.10}\\
i(\ga_j, 2\hat{m}_j)&=& 2\varrho_j \hat{p} N -(S^+_{M_j}(1)+C(M_j)-2\hat{\Delta}_j),     \lb{3.11}\eea
where $\varrho_i$ and $Q_j(m)$ are the same as those in (\ref{3.5}) and \bea \hat{m}_j=\left(\left[\frac{\hat{p}N}{M|\hat i(\gamma_j, 1)|}\right]+\hat{\chi}_j\right)M,\quad
\hat{\Delta}_j = \sum_{0<\{\hat{m}_j\th/\pi\}<\delta}S^-_{M_j}(e^{\sqrt{-1}\th}),\quad\forall\  1\le j\le q.\lb{3.12}\eea
Furthermore, comparing with some corresponding integers in Theorem 3.2, there holds
\bea \hat{\chi}_j=\chi_j,\quad \hat{m}_j=\hat{p}m_j,\quad \hat{\Delta}_j=\Delta_j,\quad \forall\ 1\le j\le q.\lb{3.13}\eea}

The following is the relation between the Viterbo index and the Maslov-type index of symplectic path.

\medskip

{\bf Theorem 3.4.} (cf. Theorem 2.1 of \cite{HuL}) {\it Suppose $\Sg\in \H_{st}(2n)$ and
$(\tau,y)\in \T(\Sigma)$. Then we have
$$ i(y^m) = i(y, m)-n,\quad \nu(y^m) =\nu(y, m),        \qquad \forall m\in\N, $$
where $i(y^m)=i(m\tau,y^m)$ and $\nu(y^m)=\nu(m\tau,y^m)$ are the index and nullity of $(m\tau,y)$ defined in (\ref{2.6})
via the Viterbo indices, $i(y, m)$ and $\nu(y, m)$ are the Maslov-type index and nullity of $(m\tau,y)$ (cf. Section 5.4
of \cite{Lon2}).}

In particular, we have $\hat{i}(\tau,y)=\hat{i}(y,1)$, where $\hat{i}(\tau ,y)$ is given in (\ref{2.7}),
$\hat{i}(y,1)$ is the mean Maslov-type index given at the end of Theorem 3.1 (cf. Definition 8.1 of \cite{Lon3}). Hence
we denote them simply by the same $\hat{i}(y)$ below.

\setcounter{figure}{0}
\setcounter{equation}{0}
\section{Proof of the main theorem}

In this section, we fix a non-degenerate $\Sg\in \H_{st}(6)$. Then $\;^{\#}\T(\Sg)\ge 2$ holds by \cite{HuL} of X. Hu and
Y. Long as mentioned before. Next we carry out two steps to give the proof of Theorem 1.1.

\medskip

{\bf Step A:} {\it The aim in this Step A is to prove the following theorem.}

{\bf Theorem 4.1.} {\it Let $\Sigma$ be a non-degenerate $C^3$ compact star-shaped hypersurface
in $\R^{6}$ possessing precisely $2$ prime closed characteristics $(\tau_1,y_1)$ and $(\tau_2,y_2)$.
Then the conclusions in the following Lemmas 4.2 to 4.7 on them hold.}

\medskip

Note that by the non-degeneracy of $\Sg$, here both $y_1$ and $y_2$ and their iterates are all non-degenerate.

Denote by $\ga_k\equiv \ga_{y_k}$, the associated symplectic path of $(\tau_k,y_k)$ for $1\le k\le 2$.
Then by Lemma 3.3 of \cite{HuL} and Theorem 3.1, there exists a path $f_k\in C([0,1],\Omega^0(\ga_k(\tau_k))$
such that $f_k(0)=\ga_k(\tau_k)$, $U_k\in Sp(4)$ such that
\be f_k(1) = N_1(1,1)\dm U_k, \qquad \forall\; 1\le k\le 2.\lb{4.1}\ee
By Theorem 3.1 and the above assumption on the non-degeneracy, we have
\bea  U_k
&=& R(\th_1)\,\dm\,\cdots\,\dm\,R(\th_r)\,\dm\,M_s\dm\,N_2(e^{\aa_{1}\sqrt{-1}},A_{1})\,\dm\,\cdots\,
                \dm\,N_2(e^{\aa_{r_{\ast}}\sqrt{-1}},A_{r_{\ast}})   \nn\\
&& \,\dm\,N_2(e^{\bb_{1}\sqrt{-1}},B_{1})\,\dm\,\cdots\,\dm\,N_2(e^{\bb_{r_{0}}\sqrt{-1}},B_{r_{0}}), \lb{4.2}\eea
where we have $\frac{\th_j}{2\pi}\not\in \Q$ for $1\le j\le r$, $M_s=D(2)^{\diamond s}$ with an integer $s\ge 0$
or $M_s=D(-2)\diamond D(2)^{\diamond (s-1)}$ with an integer $s\ge 1$, $\frac{\aa_j}{2\pi}\not\in \Q$ with
$N_2(e^{\aa_j\sqrt{-1}},A_j)$ being non-trivial for $1\le j\le r_{\ast}$, $\frac{\bb_j}{2\pi}\not\in \Q$ with
$N_2(e^{\bb_j\sqrt{-1}},B_j)$ being trivial for $1\le j\le r_0$ as defined in Definition 1.8.11 of \cite{Lon3} and
\be r+ s +2r_{\ast} + 2r_0 = n-1=2. \lb{4.3}\ee

Here when $U_k = N_2(e^{\eta\sqrt{-1}},C)$ holds for some $\frac{\eta}{2\pi}\not\in \Q$ and $2\times 2$ matrix $C$,
we call the corresponding prime closed characteristic $(\tau_k,y_k)$ being {\bf $N_2$-elliptic}.

At first we prove in the following Lemmas 4.2 to 4.5 that both $(\tau_1, y_1)$ and
$(\tau_2, y_2)$ are non-hyperbolic and non-$N_2$-elliptic. Thus we may assume that $(\tau_2, y_2)$ is hyperbolic
or $N_2$-elliptic, and want to derive a contradiction.

In particular, since $(\tau_2, y_2)$ is hyperbolic or $N_2$-elliptic by assumption, we have
\be U_2 = D(2)^{\dm 2}, \quad D(-2)\dm D(2) \quad\mbox{or}\quad N_2(e^{\eta\sqrt{-1}},C), \lb{4.4}\ee
where we have $\frac{\eta}{2\pi}\not\in \Q$ and $N_2(e^{\eta\sqrt{-1}},C)$ being non-trivial or trivial as in (\ref{4.2}).

Note that by (8.2.15) of Theorem 8.2.3 and (8.2.27) of Theorem 8.2.4 in \cite{Lon3}, we have
$\varphi(\frac{m\eta}{2\pi}) = E(\frac{m\eta}{2\pi}) - [\frac{m\eta}{2\pi}] =1$, and the Maslov-type index iteration
formula for symplectic matrix paths ending at $N_2(e^{\eta\sqrt{-1}},C)$ coincides with that of the truly hyperbolic
symplectic matrix paths in (8.2.3) of Theorem 8.2.1 of \cite{Lon3}, because $\frac{\eta}{2\pi}\not\in \Q$. Note that
because (8.2.16) in \cite{Lon3} holds in Theorems 8.2.3 and 8.2.4 in \cite{Lon3}, the initial Maslov-type index of
such a $4\times 4$ symplectic matrix path must be even. Hence by (\ref{4.1})-(\ref{4.4}), Theorems 3.1 and 3.2, we
obtain
\be \left\{\matrix{
i(y_1^m) = m(i(y_1)+n+1-r)+2\sum_{j=1}^r \left[\frac{m\theta_j}{2\pi}\right]+r-1-n, \qquad\quad \nn\cr
\qquad\qquad = m(i(y_1)+4-r)+2\sum_{j=1}^r \left[\frac{m\theta_j}{2\pi}\right]+r-4,\quad\forall\ m\ge 1, \qquad\quad \nn\cr
i(y_2^m) = m(i(y_2)+n+1)-1-n=m(i(y_2)+4)-4,\quad\forall\ m\ge 1,  \cr}\right. \lb{4.5}\ee
where in (\ref{4.5}), we have used $E(a)=[a]+1$ for $a\not\in \Z$.
Thus we have
\bea
\hat{i}(y_1) &=& i(y_1)+4-r+\sum_{j=1}^r\frac{\theta_j}{\pi}, \lb{4.6}\\
\hat{i}(y_2) &=& i(y_2)+4.  \lb{4.7}\eea

\medskip

{\bf Lemma 4.2.} {\it We have $\hat i(y_1)>0$ and $\hat i(y_2)>0$, and consequently $i(y_1)>-6$ and $i(y_2)>-4$ hold
by (\ref{4.6}) and (\ref{4.7}) respectively.}

\medskip

{\bf Proof.} By (\ref{2.9}) of Theorem 2.3, we have
\be \sum_{1\le j\le 2 \atop \hat{i}(y_j)>0}\frac{\hat{\chi}(y_j)}{\hat{i}(y_j)}=\frac{1}{2},\qquad
    \sum_{1\le j\le 2 \atop \hat{i}(y_j)<0}\frac{\hat{\chi}(y_j)}{\hat{i}(y_j)}=0,  \lb{4.8}\ee
Note that the first identity in (\ref{4.8}) implies that at least one of $\hat{i}(y_1)$ and $\hat{i}(y_2)$ is positive.
Then if the other one has negative mean index, then its Euler characteristic must be non-zero by (\ref{2.10}), which
contradicts to the second identity of (\ref{4.8}). Thus there is no closed characteristic possessing negative mean index
in our case.

Next we prove the lemma by contradiction and assume that one of the two prime closed characteristics possesses zero mean
index. We carry out the proof in two cases.

{\bf Case 1.} {\it $\hat{i}(y_1)=0$ and $\hat{i}(y_2)>0$ hold.}

By (\ref{4.6}) we have
\be i(y_1)+4-r+\sum_{j=1}^r\frac{\theta_j}{\pi}=0.  \lb{4.9}\ee
Noticing that $\frac{\theta_j}{\pi}\notin \Q$, it yields $r=0$ or 2. If $r=0$, then $i(y_1)=-4$ and together with (\ref{4.5})
it gives $i(y_1^m)=-4$ for all $m\in\N$. If $r=2$, noticing that $\sum_{j=1}^2\frac{\theta_j}{\pi}\in (0,4)$ and
$i(y_1)\in 2\Z$ by (8.1.8) of Theorem 8.1.4 and (8.1.29) of Theorem 8.1.7 of \cite{Lon3}, then by  (\ref{4.9}) we have
$\sum_{j=1}^2\frac{\theta_j}{\pi}=2$ and $i(y_1)=-4$. Since $\sum_{j=1}^2\frac{m\theta_j}{2\pi}=m$ implies that
$\sum_{j=1}^2[\frac{m\theta_j}{2\pi}]=m-1$, so by (\ref{4.5}) we have
$$ i(y_1^m)=-2m+2\sum_{j=1}^r\left[\frac{m\theta_j}{2\pi}\right]-2=-2m+2(m-1)-2=-4,\quad \forall\ m\ge 1. $$
In summary, we always have
\be i(y_1^m)=-4, \quad \forall\ m\ge 1.  \lb{4.10}\ee

Now based on (\ref{4.10}), the proof of Lemma 4.2 follows from an argument in the proof of Theorem 1.6 of \cite{DLLW2}.
For readers' conveniences we describe it briefly in six steps as follows.

{\bf Step 1.} On one hand, by (\ref{4.10}), there always holds $i(y_1^m)=-4$ for any $m\in\N$. On the
other hand, note that $\hat{i}(y_2)>0$ implies $i(y_2^m)\to +\infty$ as
$m\to +\infty$. Thus iterates $y_2^m$ have indices satisfying $i(y_2^m)\neq -3$, $-4$ and
$-5$ for any large enough $m\in\N$. Therefore for large enough $a$, all the closed characteristics $y_j^m$ for
$1\le j\le 2$ with period larger than $aT$, which implies that the iterate number $m$ is very large, will have their
Viterbo indices:
\be \left\{\matrix{
& \mbox{either (i) equal to} -4, \quad \mbox{when}\; \hat{i}(y_j)=0,  \cr
& \mbox{or (ii) different from} -3,\ -4\ \mbox{and} -5, \quad \mbox{when}\; \hat{i}(y_j)\not= 0.\cr}\right. \lb{4.11}\ee

\medskip

{\bf Step 2.} For $a\in\R$, let $X^-(a,K)=\{x\in X\mid F_{a,K}(x)<0\}$ with $K=K(a)$ as defined in Section 2
as well as in Section 7 of \cite{Vit2}. Note also that the origin $0$ of $X$ is not contained in $X^-(a,K)$ by definition.
Because the Hamiltonian function $H_{a,K}$ is quadratic homogeneous as assumed at the beginning of Section 7 of
\cite{Vit2} due to the study there being near the origin, the functional $F_{a,K}$ is homogeneous too.

For any large enough positive $a<a'$, we fix the same constant $K'>0$ as that in the second Step on p.639 of\cite{Vit2}
to be sufficiently large than $K$ such that the Hamiltonian function $H_{t,K'}(x)$ is strictly convex for every
$t\in [a,a']$. Now let $A=X^-(a,K')$ and $A'=X^-(a',K')$.
Because the period set $P_{\Sg}= \{m\tau_j\;|\;1\le j\le 2, m\in\N\}$  is discrete, we choose the above constants $a$
and $a'$ carefully such that $aT$ and $a'T$ do not belong to $P_{\Sg}$. Note that for $t\in [a,a']$ because every
critical orbit $S^1\cdot x$ of the functional $F_{t,K'}$ always possesses the critical value $F_{t,K'}(S^1\cdot x) = 0$ as
mentioned in p.639 of \cite{Vit2} and by (2.7) of \cite{LLW}, the boundary sets of $A$ and $A'$, i.e.,
$\{x\,|\,F_{t,K'}(x)=0\}$ with $t=a$ or $a'$, possess no critical orbits, and specially the origin $0$ of $X$ is not
contained in $A$ and $A'$. Therefore by the homogeneity mentioned above we have
\bea   H_{S^1, d(K')+i}(A',A) \;=\; H_{S^1, d(K')+i}(A'\cap S(X),A\cap S(X)), \qquad \forall\;i\in \Z, \lb{4.12}\eea
where $S(X)$ is the unit sphere of $X$. Note that a closed
characteristic whose period locates between $aT$ and $a'T$ corresponds to a
critical orbit of $F_{t,K'}$ for some $t\in [a,a']$ contained in $(A'\bs A)\cap S(X)$.

\medskip

{\bf Step 3.} Now for the chosen large enough $a$ and $a'$ with $a<a'$, by (\ref{4.11}) there exists no any closed
characteristic whose period locates between $aT$ and $a'T$ possessing Viterbo index $-3$ or $-5$. Therefore by the
argument in Step 3 in the proof of Theorem 1.6 of \cite{DLLW2}, we obtain
\be  H_{S^1, d(K')-3}(A',A) \;=\; H_{S^1, d(K')-5}(A',A) \;=\; 0.   \lb{4.13}\ee

\medskip

{\bf Step 4.} Now we consider the following exact sequence of the triple $(X,A',A)$
\bea
&&\cdots\longrightarrow H_{S^1, d(K')-3}(A',A)\stackrel{i_{3*}}{\longrightarrow} H_{S^1, d(K')-3}(X,A)
     \stackrel{j_{3*}}{\longrightarrow}H_{S^1, d(K')-3}(X,A')\nn\\
&&\quad\qquad\stackrel{\partial_{3*}}{\longrightarrow} H_{S^1, d(K')-4}(A',A)\stackrel{i_{4*}}{\longrightarrow} H_{S^1, d(K')-4}(X,A)
                         \stackrel{j_{4*}}{\longrightarrow}H_{S^1, d(K')-4}(X,A')\nn\\
&&\qquad\qquad\stackrel{\partial_{4*}}{\longrightarrow} H_{S^1, d(K')-5}(A',A)\longrightarrow\cdots. \lb{4.14}\eea

Then we consider the following homomorphisms:
\bea
H_{S^1, d(K)-3}(X,X^-(a,K)) &{\xi_1}\atop{\longrightarrow}& H_{S^1, d(K')-3}(X,X^-(a,K')),  \nn\\
H_{S^1, d(K')-3}(X,X^-(a,K')) &{\xi_2}\atop{\longrightarrow}& H_{S^1, d(K')-3}(X,X^-(a',K')), \nn\\
H_{S^1, d(K)-3}(X,X^-(a,K)) &{\xi}\atop{\longrightarrow}& H_{S^1, d(K')-3}(X,X^-(a',K')), \nn\eea
where $\xi_1$ is the homomorphism given by (7.2) of \cite{Vit2}, $j_{3*}=\xi_2$ is the homomorphism
given by the line above (7.4) of \cite{Vit2}, and $\xi= \xi_2\circ\xi_1$ is precisely the homomorphism
given by (7.4) of \cite{Vit2}. Here $\xi_1$ is an isomorphism and $\xi$ is a zero homomorphism as
proved in the Steps 1 and 2 of the proof of Theorem 7.1 in \cite{Vit2} respectively. Therefore
$j_{3*}=\xi_2$ is also an zero homomorphism.

Therefore (\ref{4.13}) and (\ref{4.14}) yield
\be  H_{S^1, d(K')-3}(X,A) = \mbox{Ker}(j_{3*}) = \mbox{Im}(i_{3*}) = i_{3*}(H_{S^1, d(K')-3}(A',A)) = 0.  \lb{4.15}\ee

Now we fix the above chosen $a'>0$ and choose another large enough $a''>a'$, and enlarge the constant $K'$
chosen above (\ref{4.12}) so that the conclusions between (\ref{4.11}) and (\ref{4.12}) hold when we replace $(a,a')$ by
$(a',a'')$. Then repeating the above proof with the long exact sequence of the triple
$(X,A'',A')$ instead of $(X,A',A)$ in the above arguments with $A''=X^-(a'',K')$, similarly we obtain
\be  H_{S^1, d(K')-3}(X,A') =0.  \lb{4.16}\ee

Together with (\ref{4.15}) and (\ref{4.16}), (\ref{4.14}) yields
\be 0\stackrel{\partial_{3*}}{\longrightarrow} H_{S^1, d(K')-4}(A',A)\stackrel{i_{4*}}{\longrightarrow}
   H_{S^1, d(K')-4}(X,A)\stackrel{j_{4*}}{\longrightarrow}H_{S^1, d(K')-4}(X,A')
    \stackrel{\partial_{4*}}{\longrightarrow} 0.  \lb{4.17}\ee

\medskip

{\bf Step 5.} When $a$ increases, we always meet infinitely many closed characteristics with Viterbo index $-4$ due to
(\ref{4.10}). For the above chosen large enough $a < a'$, there exist only finitely many closed characteristics among
$\{y_1^m\ |\  m\ge 1\}$ such that their periods locate between $aT$ and $a'T$. As Step 5 in the proof of Theorem 1.6
of \cite{DLLW2}, we obtain
\bea  H_{S^1, d(K')-4}(A',A) \not= 0, \lb{4.18}\eea

\medskip

{\bf Step 6.} By the exactness of the sequence (\ref{4.17}) and (\ref{4.18}), we obtain
$$   H_{S^1, d(K')-4}(X,A)=H_{S^1, d(K')-4}(A',A)\bigoplus H_{S^1, d(K')-4}(X,A')\neq 0.  $$
Then, by our choice of $a$, $a'$ and $a''$, and replacing $(X,A',A)$ by $(X,A'',A')$ in the above arguments, similarly
we obtain
\be H_{S^1, d(K')-4}(X,A')\neq 0. \lb{4.19}\ee

Now on one hand, if $j_{4*}$ in (\ref{4.17}) is a trivial homomorphism, then by the exactness of the sequence
(\ref{4.17}) it yields
$$   H_{S^1, d(K')-4}(X,A') = \mbox{Ker}(\partial_{4*}) = \mbox{Im}(j_{4*})=0,   $$
which contradicts to (\ref{4.19}). Therefore $j_{4*}$ in (\ref{4.17}) is a non-trivial homomorphism.

However, on the other hand, by our discussion between (\ref{4.14}) and (\ref{4.15}) using the argument in \cite{Vit2},
$j_{4*}$ in (\ref{4.17}) is a zero homomorphism. This contradiction completes the proof of case (i) of Lemma 4.2.

{\bf Case 2.} {\it $\hat{i}(y_2)=0$ and $\hat{i}(y_1)>0$ hold.}

Then by (\ref{4.7}) we have $i(y_2)=-4$, together with (\ref{4.5}) which gives $i(y_2^m)=-4$ for all $m\in\N$, then
we get a contradiction too as in the above Case 1.

Note that the last claim on $i(y_1)$ and $i(y_2)$ follows from the above results $\hat{i}(y_1)>0$ and $\hat{i}(y_2)>0$
and (\ref{4.5})-(\ref{4.7}), and Lemma 4.2 is proved. \hfill\hb

\medskip

Now assuming $^{\#}\T(\Sg)=2$, we prove that both the two prime closed characteristics $(\tau_1,y_1)$ and $(\tau_2,y_2)$
are non-hyperbolic and non-$N_2$-elliptic through Lemmas 4.3-4.5 below according to the numbers of irrational rotations
appeared in the decomposition of $U_1$ in (\ref{4.2}), and then obtain further properties about $y_1$ and $y_2$ through
Lemmas 4.6-4.7.

\medskip

{\bf Lemma 4.3.} {\it If $r=0$ in (\ref{4.3}) for $U_1$, we can derive a contradiction.}

\medskip

{\bf Proof.} In this case, we have $i(y_1^m)=m(i(y_1)+4)-4$ by (\ref{4.5}).

Note that here $y_1$ can be hyperbolic or $N_2$-elliptic too. Then together with (\ref{2.10}) we have $\hat{\chi}(y_j)=1$
if $i(y_j)\in 2\Z$, and $\hat{\chi}(y_j)=-\frac{1}{2}$ if $i(y_j)\in 2\Z+1$ for $j=1, 2$. Then in order to make the first
identity in (\ref{4.8}) with $k=2$ hold, $i(y_1)=0=i(y_2)$ must hold. Thus we have $i(y_j^m)=4m-4\ge 4$ for $m\ge 2$ and
$j=1,2$ by (\ref{4.5}). Then $M_2=0$ by (\ref{2.8}) and (\ref{2.11}). Together with (\ref{2.12}) it yields a
contradiction $0=M_2\ge b_2=1$, and completes the proof of the Lemma 4.3. \hfill\hb

\medskip

Note that Lemma 4.3 shows also that the case in which both $(\tau_1,y_1)$  and $(\tau_2,y_2)$ are hyperbolic or $N_2$-elliptic
can not happen.

\medskip

{\bf Lemma 4.4.} {\it If $r=1$ in (\ref{4.3}) for $U_1$, we can derive a contradiction.}

\medskip

{\bf Proof.} Note that $\hat{\chi}(y_i)\in\Q$ for $i=1$ and $2$ in this case. Since $r=1$, we have $\hat{i}(y_1)\notin\Q$. Thus
by the first identity in (\ref{4.8}) with $k=2$, we obtain $\hat{i}(y_2)\notin\Q$. But this contradicts to
$\hat{i}(y_2)=i(y_2)+4 \in \Q$. \hfill\hb

\medskip

{\bf Lemma 4.5.} {\it If $r=2$ in (\ref{4.3}) for $U_1$, we can derive a contradiction.}

\medskip

{\bf Proof.} By Theorem 2.3, we have
\be \frac{1}{2} = \frac{\hat{\chi}(y_1)}{\hat{i}(y_1)}+\frac{\hat{\chi}(y_2)}{\hat{i}(y_2)}
     = \frac{1}{i(y_1)+2+\sum_{j=1}^2\frac{\theta_j}{\pi}}+\frac{\hat{\chi}(y_2)}{\hat{i}(y_2)}. \lb{4.20}\ee
By Lemma 4.2, We have $i(y_1)>-6$ and $i(y_2)\ge -3$. By Theorems 3.1 and 3.2, we have $i(y_1)\in 2\Z$. Then
$i(y_1)\ge -4$ holds. We continue our study below in three cases according to the value of $i(y_1)$.

\medskip

{\bf Case 1.} {\it $i(y_1)=-4$.}

\medskip

In this case, if $i(y_2)\ge -2$, then we have $M_{-4}\ge 1$ and $M_{-3}=0$ by Lemma 2.2. Thus we have
\bea -1
&\ge& M_{i(y_1)+1}-M_{i(y_1)}+M_{i(y_1)-1}-M_{i(y_1)-2} \nn\\
&\ge& b_{i(y_1)+1}-b_{i(y_1)}+b_{i(y_1)-1}-b_{i(y_1)-2} = 0-0+0-0 = 0. \lb{4.21}\eea
This contradiction yields $i(y_2)=-3$ and then
\bea
i(y_1^m) &=& -2m+2\sum_{j=1}^r\left[\frac{m\theta_j}{2\pi}\right]-2,\qquad\forall\ m\ge 1,\lb{4.22}\\
i(y_2^m) &=& m-4,\qquad\forall\ m\ge 1.\lb{4.23}\eea
Therefore (\ref{4.20}) becomes
$$  \frac{1}{-2+\sum_{j=1}^2\frac{\theta_j}{\pi}}-\frac{1}{2}=\frac{1}{2}. $$
Hence
\be \sum_{j=1}^2\frac{\theta_j}{\pi}=3. \lb{4.24}\ee
Denote by $\alpha_j=\frac{\theta_j}{2\pi}\in(0,\,1)\setminus\Q$. Then by (\ref{4.22}) and (\ref{4.24}) we have
\bea i(y_1^m)
&=& -2m+2\sum_{j=1}^2\left[\frac{m\theta_j}{2\pi}\right]-2   \nn\\
&=& -2m+2\sum_{j=1}^2 [m\alpha_j]-2    \nn\\
&=& -2m+2[m\alpha_1]+2\left[m\left(\frac{3}{2}-\alpha_1\right)\right]-2  \nn\\
&=& -2m+ 2[m\alpha_1]+2\left[m+\frac{m}{2}-([m\alpha_1]+\{m\alpha_1\})\right]-2  \nn\\
&=& 2\left[\frac{m}{2}-\{m\alpha_1\}\right]-2,\quad\forall\ m\ge 1.   \lb{4.25}\eea

{\bf Claim 1.} {\it We have
\bea
M_{2k+1} &=& \left\{\matrix{1, &\quad {\rm if}\;\; 2k+1\ge -3, \cr
                            0, &\quad {\rm if}\;\; 2k+1<-3. \cr}\right.\lb{4.26}\\
M_{2k} &=& \left\{\matrix{1, &\quad {\rm if}\;\; -4\le 2k<0, \cr
                          2, &\quad {\rm if}\;\; 2k\ge 0, \cr
                          0, &\quad {\rm if}\;\; 2k<-4, \cr}\right.\lb{4.27}\\
i(y_1^{2m+1}) &=& 2m-2,\qquad \forall\ m\ge 1. \lb{4.28}\eea}

In fact, (\ref{4.26}) follows from (\ref{4.23}) and Lemma 2.2. By (\ref{4.25}), we have
\bea
i(y_1^{2m}) &=& 2m-4,\qquad\forall\ m\ge 1,\lb{4.29}\\
i(y_1^{2m+1}) &=& 2m-2+2\left[\frac{1}{2}-\{(2m+1)\alpha_1\}\right],\qquad\forall\ m\ge 1.\lb{4.30}\eea
Hence
\be 2m-4\le i(y_1^{2m+1})\le 2m-2,\qquad\forall\ m\ge 1. \lb{4.31}\ee
Thus $i(y_1^3)=-2$ or $i(y_1^3)=0$ holds. If $i(y_1^3)=-2$, then $M_{-2}\ge 2$ by Lemma 2.2 and (\ref{4.29}).
And then by (\ref{2.12}) it yields
\bea -1
&=& 1-2+1-1\ge M_{-1}-M_{-2}+M_{-3}-M_{-4} \nn\\
&\ge& b_{-1}-b_{-2}+b_{-3}-b_{-4}=0. \nn\eea
Thus $i(y_1^3)=0$ must hold and then
\be M_0\ge 2,\qquad M_{-4}=M_{-2}=1. \lb{4.32}\ee

By (\ref{4.31}), $i(y_1^5)=0$ or $i(y_1^5)=2$. If $i(y_1^5)=0$, then $M_{0}\ge 3$ by (\ref{4.32}).
Then by (\ref{2.12})
\bea -2
&=& 1-3+1-1+1-1 \ge M_1-M_0+M_{-1}-M_{-2}+M_{-3}-M_{-4} \nn\\
&\ge& b_1-b_0+b_{-1}-b_{-2}+b_{-3}-b_{-4}=-1. \nn\eea
Thus $i(y_1^5)=2$ must hold and then $M_0=2$ and $M_2\ge 2$.

Now we use an induction argument on $k$. For $k\ge 3$, we assume
\be i(y_1^{2m+1})=2m-2, \forall\ 1\le m\le k-1,\quad M_{2j}=2,\ \forall \ 0\le j\le k-3,\quad M_{2k-4}\ge 2.\lb{4.33}\ee

Then by (\ref{4.31}), there holds $i(y_1^{2k+1})=2k-2$ or $i(y_1^{2k+1})=2k-4$.
If $i(y_1^{2k+1})=2k-4$, then by the induction assumption (\ref{4.33}) and (\ref{2.12}) we have
\bea -k
&=& 1-3+\sum_{j=0}^{k-3}(1-2)+1-1+1-1  \nn\\
&\ge& M_{2k-3}-M_{2k-4}+\sum_{j=0}^{k-3}(M_{2j+1}-M_{2j})+M_{-1}-M_{-2}+M_{-3}-M_{-4}  \nn\\
&\ge& \sum_{j=-2}^{k-2}(b_{2j+1}-b_{2j})=-\sum_{j=0}^{k-2} b_{2j}=-k+1.\nn\eea
This contradiction proves $i(y_1^{2k+1})=2k-2$. And then $M_{2k-4}=2$, (\ref{4.27}) and (\ref{4.28}) hold.
Hence Claim 1 holds.

\medskip

On the other hand, by (\ref{4.30}) and the Kronecker's uniform distribution theorem about irrational
numbers, there must exist infinitely many $m\in\N$ such that $i(y_1^{2m+1})=2m-4$. This contradicts
to $i(y_1^{2m+1})=2m-2$ in Claim 1, which shows that Case 1 can not happen.

\medskip

{\bf Case 2.} {\it $i(y_1)=-2$.}

\medskip

In this case, if $i(y_2)\ge 0$, then we have $M_{-2}\ge 1$ and $M_{-1}=0$ by Lemma 2.2. Thus we still have
(\ref{4.21}). This contradiction yields $i(y_2)\in\{-3, -2, -1\}$ and then
\bea
i(y_1^m) &=& 2\sum_{j=1}^r \left[\frac{m\theta_j}{2\pi}\right]-2,\quad\forall\ m\ge 1, \lb{4.34}\\
i(y_2^m) &=& m(i(y_2)+4)-4,\quad\forall\ m\ge 1.  \lb{4.35}\eea
We continue the proof in three subcases according to the value of $i(y_2)$.

{\bf Subcase 2.1.} {\it $i(y_2)=-3$.}

Therefore (\ref{4.20}) becomes
$$ \frac{1}{\sum_{j=1}^2\frac{\theta_j}{\pi}}-\frac{1}{2}=\frac{1}{2}. $$
Hence
\be \sum_{j=1}^2\frac{\theta_j}{\pi}=1.  \lb{4.36}\ee

Denote by $\alpha_j=\frac{\theta_j}{2\pi}\in(0,\,1)\setminus\Q$. Then by (\ref{4.34}) and (\ref{4.36})
we have
\bea  i(y_1^m)\
&=& 2\sum_{j=1}^2 \left[\frac{m\theta_j}{2\pi}\right]-2  \nn\\
&=& 2\sum_{j=1}^2 [m\alpha_j]-2   \nn\\
&=& 2[m\alpha_1]+2\left[m\left(\frac{1}{2}-\alpha_1\right)\right]-2  \nn\\
&=& 2\left[\frac{m}{2}-\{m\alpha_1\}\right]-2,\quad\forall\ m\ge 1.  \lb{4.37}\eea

{\bf Claim 2.} {\it We have
\bea
M_{2k+1} &=& \left\{\matrix{1, &\qquad {\rm if}\;\; 2k+1\ge -3, \cr
                            0, &\qquad {\rm if}\;\; 2k+1<-3. \cr}\right.\lb{4.38}\\
M_{2k} &=& \left\{\matrix{2, &\qquad {\rm if}\;\; k\ge -1, \cr
                          0, &\qquad {\rm if}\;\; k<-1, \cr}\right.\lb{4.39}\\
i(y_1^{2m+1}) &=& 2m-2,\qquad \forall\ m\ge 1. \lb{4.40}\eea}

In fact, (\ref{4.38}) follows from (\ref{4.35}) and Lemma 2.2.
By (\ref{4.37}), we have
\bea
i(y_1^{2m}) &=& 2m-4,\quad \forall \  m\ge 1,  \lb{4.41}\\
i(y_1^{2m+1}) &=& 2m-2+2\left[\frac{1}{2}-\{(2m+1)\alpha_1\}\right],\quad\forall\ m\ge 1. \lb{4.42}\eea
Hence we obtain
$$ 2m-4\le i(y_1^{2m+1})\le 2m-2,\qquad\forall\ m\ge 1. $$
Thus $i(y_1^3)=-2$ or $i(y_1^3)=0$ holds. If $i(y_1^3)=-2$, then $M_{-2}\ge 3$ holds by Lemma 2.2 and
(\ref{4.41}). And then
\bea -1
&=& 1-3+1 \ge M_{-1}-M_{-2}+M_{-3}  \nn\\
&\ge& b_{-1}-b_{-2}+b_{-3}=0.  \nn\eea
Thus $i(y_1^3)=0$ must holds and then we obtain
$$  M_0\ge 2,\qquad M_{-2}=2.  $$

Then by a similar argument as that in the proof of Claim 1, we obtain (\ref{4.39}) and (\ref{4.40}),
and complete the proof of Claim 2.

On the other hand, by (\ref{4.42}) and the Kronecker's uniform distribution theorem, there must exist
infinitely many $m\in\N$ such that $i(y_1^{2m+1})=2m-4$. This contradicts to the fact $i(y_1^{2m+1})=2m-2$
in Claim 2. Thus the Subcase 2.1 can not happen.

\medskip

{\bf Subcase 2.2.} {\it $i(y_2)=-2$.}

\medskip

Then (\ref{4.20}) becomes
$$   \frac{1}{\sum_{j=1}^2\frac{\theta_j}{\pi}}+\frac{1}{2}=\frac{1}{2}. $$
Hence we get a contradiction, which shows that the Subcase 2.2 can not happen.

\medskip

{\bf Subcase 2.3.} {\it $i(y_2)=-1$.}

\medskip

Therefore (\ref{4.20}) becomes
$$ \frac{1}{\sum_{j=1}^2\frac{\theta_j}{\pi}}-\frac{1}{6}=\frac{1}{2}. $$
Hence we obtain
\be  \sum_{j=1}^2\frac{\theta_j}{\pi}=\frac{3}{2}.  \lb{4.43}\ee

Denote by $\alpha_j=\frac{\theta_j}{2\pi}\in(0,\,1)\setminus\Q$. Then by (\ref{4.34}) and (\ref{4.43})
we have
\bea i(y_1^m)
&=& 2\sum_{j=1}^2 \left[\frac{m\theta_j}{2\pi}\right]-2  \nn\\
&=& 2\sum_{j=1}^2 [m\alpha_j]-2  \nn\\
&=& 2[m\alpha_1]+2\left[m\left(\frac{3}{4}-\alpha_1\right)\right]-2  \nn\\
&=& 2\left[\frac{3m}{4}-\{m\alpha_1\}\right]-2,\quad\forall\ m\ge 1.  \lb{4.44}\eea
Then by (\ref{4.35}) and Lemma 2.2 we obtain
\be M_{2k-1} = \left\{\matrix{1, &\quad {\rm if}\;\; k\in3\N_0, \cr
                              0, &\quad {\rm otherwise.}\cr}\right. \lb{4.45}\ee
By (\ref{4.44}), we have
\bea
i(y_1^{4m}) &=& 6m-4,\qquad\forall\ m\ge 1,\nn\\
i(y_1^{4m+k}) &=& 6m-2+2\left[\frac{3k}{4}-\{(4m+k)\alpha_1\}\right],\quad \forall \ 1\le k\le 3,\  m\ge 1. \nn\eea
Hence we obtain
\be   6m-4+2(k-1)\le i(y_1^{4m+k})\le 6m-2+2(k-1),\quad \forall \ 1\le k\le 3,\  m\ge 1.\lb{4.46}\ee

Thus $i(y_1^2)=-2$ or $i(y_1^2)=0$ holds. If $i(y_1^2)=-2$, then $M_{-2}\ge 2$ by $i(y_1)=-2$ and Lemma 2.2,
and then
$$ -1=1-2\ge M_{-1}-M_{-2}\ge b_{-1}-b_{-2}=0.  $$
Thus $i(y_1^2)=0$ must hold and then
\be M_0\ge 1,\qquad M_{-2}=1.  \lb{4.47}\ee

Note that $i(y_1^3)=0$ or $i(y_1^3)=2$ holds by (\ref{4.46}) again. If $i(y_1^3)=0$, then $M_{0}\ge 2$ by
Lemma 2.2 and (\ref{4.47}). And then we have
$$ -2=0-2+1-1\ge M_1-M_0+M_{-1}-M_{-2}\ge b_1-b_0+b_{-1}-b_{-2}=-1. $$
Thus $i(y_1^3)=2$ must hold and then we obtain
\be M_2\ge 2,\qquad M_{0}=1.  \lb{4.48}\ee

Therefore by (\ref{4.45}) and (\ref{4.48}) it yields
\bea -3 = 0-2+0-1+1-1
&\ge& M_3-M_2+M_1-M_0+M_{-1}-M_{-2} \nn\\
&\ge& b_3-b_2+b_1-b_0+b_{-1}-b_{-2}=-2.  \nn\eea
This contradiction shows that the Subcase 2.3 can not happen.

Therefore the Case 2 can not happen.

\medskip

{\bf Case 3.} {\it $i(y_1)\ge 0$.}

\medskip

In this case, we have
$$  \hat{i}(y_1)=i(y_1)+2+\sum_{j=1}^r\frac{\theta_j}{\pi}>2. $$
Together with $\hat{\chi}(y_1)=1$ it yields $\frac{\hat{\chi}(y_1)}{\hat{i}(y_1)}<\frac{1}{2}$. Hence
$\frac{\hat{\chi}(y_2)}{\hat{i}(y_2)}>0$ by (\ref{4.20}), which implies that $i(y_2)\in 2\Z$.

Then we have $i(y_2)\ge -2$ by Lemma 4.1. If $i(y_2)=-2$, then we have $M_{-2}=1$ and $M_{-1}=0$. Hence
by (\ref{2.12}) it yields
$$   -1 = M_{-1}-M_{-2} \ge b_{-1}-b_{-2} = 0.  $$
This contradiction implies $i(y_2)\ge 0$. Then by Theorem 3.2 we obtain $i(y_2,1)\ge 3$ and $i(y_1,1)\ge 3$.
So the index perfect condition of Theorem 1.2 in \cite{DLLW1} is satisfied, and then there exist at least
three prime closed characteristics. It contradicts to the assumption $^{\#}\T(\Sg)=2$, and proves that the
Case 3 can not happen.

The proof of Lemma 4.5 is complete.\hfill\hb

\medskip

Now it follows from Lemmas 4.3-4.5 that both of $(\tau_1,y_1)$ and $(\tau_2,y_2)$ are non-hyperbolic and
non-$N_2$-elliptic.

\medskip

{\bf Lemma 4.6.} {\it When $^{\#}\T(\Sg)=2$, Denote the two prime closed characteristics by $(\tau_1,y_1)$
and $(\tau_2,y_2)$. Then one of the following two cases must happen.

(i) Both $(\tau_1,y_1)$ and $(\tau_2,y_2)$ are non-hyperbolic and non-elliptic, and the matrix $U_j$ in the
decomposition (\ref{4.1}) of the end matrix $\ga_j(\tau_j)$ of the symplectic matrix path $\ga_j$ associated
to $(\tau_j,y_j)$ with $j=1$ or $2$ satisfies
\be   U_1 = R(\th_1)\dm D(2), \qquad U_2 = R(\th_2)\dm D(-2),  \lb{4.49}\ee
for some $\th_j$ satisfying $\frac{\th_j}{2\pi}\in \R\bs\Q$ with $j=1, 2$.

(ii) $(\tau_1,y_1)$ is irrationally elliptic, $(\tau_2,y_2)$ is non-hyperbolic and non-elliptic, and the
matrix $U_j$ in the decomposition (\ref{4.1}) of the end matrix $\ga_j(\tau_j)$ of the symplectic matrix
path $\ga_j$ associated to $(\tau_j,y_j)$ with $j=1$ or $2$ satisfies
\be   U_1 = R(\th_1)\dm R(\th_2), \qquad U_2 = R(\th_3)\dm D(2),  \lb{4.50}\ee
for some $\th_k\in (0,2\pi)$ satisfying $\frac{\th_k}{2\pi}\in \R\bs\Q$ with $k=1, 2, 3$.}

\medskip

{\bf Proof.} Denote the two prime closed characteristics on $\Sg$ by $(\tau_1,y_1)$ and $(\tau_2,y_2)$.
By Lemmas 4.3-4.5, both of them are non-hyperbolic and non-$N_2$-elliptic. Note that $\Sg$ is non-degenerate,
so both of the two prime closed characteristics and all of their iterates are non-degenerate. Then by the
basic normal form decomposition of symplectic matrices introduced in \cite{Lon2} (cf. Section 1.8 and
Chapter 8 of \cite{Lon3}), using notations in (\ref{4.1}) and (\ref{4.2}) we have $U_j$ with $j=1, 2$ may
be a $\dm$-sum of two of the following matrices,
$$  R(\th_1)\dm R(\th_2),\qquad R(\th_3)\dm D(a), $$
with $a\in\{2, -2\}$, $\frac{\th_k}{2\pi}\in \R\bs\Q$ with $k=1, 2, 3$.

Next we carry our the proof of Lemma 4.6 by contradiction in four cases according to the form of $U_1$
and $U_2$.

\medskip

{\bf Case 1.} {\it $U_j = R(\th_{j,1})\dm R(\th_{j,2})$ with $\frac{\th_{j,k}}{2\pi} \in \R\bs\Q$ for
$j=1, 2$ and $k=1, 2$.}

\medskip

In this case by (8.1.7)-(8.1.8) of Theorem 8.1.4 of \cite{Lon3} on the symplectic matrix path ending at
$N_1(1,1)$, and (8.1.27)-(8.1.29) of Theorem 8.1.7 of \cite{Lon3} on the symplectic matrix path ending at
$R(\th)$, each of these single paths and its iterates possesses odd Maslov-type index, and each $\dm$-sum of
their iterates of these three single paths possesses always odd Maslov-type index. Then by Theorem 3.2, we obtain
\be  i(y_j^m)\in 2\Z, \qquad \forall\;m\in\N,\ j=1,2.  \lb{4.51}\ee

Note that by (\ref{2.5}) both $d(K)$ and $d(pK)$ for every $p\in\N$ are even and then the Morse index of the
functional $F_{a,K}$ at each iterate of $y_j$ must be even. Note that in this case, $\beta(y_j^m)=1$ holds always
in (\ref{2.8}). Thus together with Lemma 2.2, the Morse type numbers defined in (\ref{2.12}) satisfy $M_{2k-1}=0$
for every $k\in\Z$. Therefore the identity in (\ref{2.12}) becomes $M(t)=\frac{1}{1-t^2}$, and then specially only
those terms of non-negative orders in (\ref{2.11}) left. Thus together with (\ref{2.5}), (\ref{2.8}) of Lemma 2.2
and (\ref{2.11}), it yields that
\be    i(y_j^m)\ge 0, \qquad \forall\;m\in\N,\;j=1,2.   \lb{4.52}\ee
Then by Theorem 3.2, the Maslov-type index of $(y_j,m)$ satisfies
\be  i(y_j,m) = i(y_j^m)+3 \ge 3, \qquad \forall\;m\in\N, \;j=1, 2.   \lb{4.53}\ee
Therefore the index perfect condition of Theorem 1.2 in \cite{DLLW1} is satisfied and it then yields
$^{\#}\T(\Sg)\ge 3$, which contradicts to the assumption $^{\#}\T(\Sg)=2$.

\medskip

{\bf Case 2.} {\it $U_j = R(\th_j)\dm D(2)$ with $\frac{\th_j}{2\pi} \in \R\bs\Q$ for $j=1, 2$.}

\medskip

In this case, by (8.2.1) and (8.2.3) of Theorem 8.2.1 in \cite{Lon3} on the symplectic matrix path ending at $D(2)$,
each of its iterates possesses always even Maslov-type index. Together with our discussions in Case 1 on symplectic
matrix paths ending at $N_1(1,1)$ and on the symplectic matrix paths ending at $R(\th)$, we obtain that each
$\dm$-sum of their iterates of these three single paths possesses always even Maslov-type index. Then by Theorem 3.2,
we obtain
\be  i(y_j^m)\in 2\Z+1, \qquad \forall\;m\in\N,\;j=1,2. \lb{4.54}\ee
Thus we can repeat our discussions in Case 1 to get that the Morse type numbers defined in (\ref{2.12}) satisfy
$M_{2k}=0$ for every $k\in\Z$. Obviously this gives a contradiction $0=M_0\ge b_0=1$ by (\ref{2.12}).

\medskip

{\bf Case 3.} {\it $U_j = R(\th_j)\dm D(-2)$ with $\frac{\th_j}{2\pi} \in \R\bs\Q$ for $j=1, 2$.}

\medskip

In this case, by (8.2.2) and (8.2.3) of Theorem 8.2.1 in \cite{Lon3} on the symplectic matrix path ending at $D(-2)$,
each of its odd iterates possesses always odd Maslov-type index, and each of its even iterates possesses always even
Maslov-type index. Together with our discussions in Case 1 on symplectic matrix paths ending at $N_1(1,1)$ and on the
symplectic matrix paths ending at $R(\th)$, we obtain that each odd iterates of the $\dm$-sum of these three single
paths possesses always odd Maslov-type index, and each even iterates of the $\dm$-sum of these three single paths
possesses always even Maslov-type index. Then by Theorem 3.2, we obtain
\be  \left\{\matrix{
    i(y_j^{2m-1})\in 2\Z, \qquad \forall\;m\in\N,\;j=1,2, \cr
    i(y_j^{2m})\in 2\Z+1, \qquad \forall\;m\in\N,\;j=1,2, \cr}\right.  \lb{4.55}\ee
Thus we obtain $\beta(y_j^{2m-1})=1$ and $\beta(y_j^{2m})=-1$ for $m\in\N$ in (\ref{2.8}). Thus every iterate of
$y_j$ contributes nothing to Morse type numbers $M_{2k-1}$, i.e., we obtain $M_{2k-1}=0$ for every $k\in\Z$. Then
we can repeat the last arguments in Case 1 to get $^{\#}\T(\Sg)\ge 3$, which contradicts to the
assumption $^{\#}\T(\Sg)=2$.

\medskip

{\bf Case 4.} {\it $U_1 = R(\th_1)\dm R(\th_2)$ and $U_2 = R(\th_3)\dm D(-2)$ with $\frac{\th_k}{2\pi} \in \R\bs\Q$
for $k=1, 2, 3$.}

\medskip

In this case, by our discussion in Case 1, we obtain (\ref{4.51}) with $j=1$ for $(\tau_1,y_1)$. Then our proof in
Case 1 below (\ref{4.51}) yields that $(\tau_1,y_1)$ contributes nothing to the Morse type number $M_{2k-1}$ for any
$k\in\Z$. By our discussion in Case 3, we obtain (\ref{4.55}) with $j=2$ for $(\tau_2,y_2)$. Then our proof in
Case 3 below (\ref{4.55}) yields that $(\tau_2,y_2)$ contributes nothing to the Morse type number $M_{2k-1}$ for any
$k\in\Z$. Summarizing the contributions of $(\tau_1,y_1)$ and $(\tau_2,y_2)$, we obtain $M_{2k-1}=0$ for any $k\in\Z$.

Then we can repeat the proof in Case 3 below (\ref{4.55}) to get $^{\#}\T(\Sg)\ge 3$, which contradicts to the
assumption $^{\#}\T(\Sg)=2$.

Thus all these four cases can not happen, and one of (\ref{4.49}) and (\ref{4.50}) must hold when $^{\#}\T(\Sg)=2$.
The proof of Lemma 4.6 is complete. \hfill\hb

\medskip

{\bf Lemma 4.7.} {\it When $^{\#}\T(\Sg)=2$, both of them must possess positive irrational mean indices.}

\medskip

{\bf Proof.} Denote the two prime closed characteristics by $(\tau_1,y_1)$ and $(\tau_2,y_2)$. Then by Theorem 1.6
of \cite{DLLW2} their mean indices satisfy $\hat{i}(y_j)\not= 0$ for $j=1, 2$. By the first identity in (\ref{4.8}),
one of them must be positive, say $\hat{i}(y_1)>0$. Because both of the two orbits and their iterates are
non-degenerate, both of the average Euler characteristics satisfy $\hat{\chi}(y_j)\not= 0$ for $j=1, 2$ by (\ref{2.10}).
Thus if $\hat{i}(y_2)<0$ held, it would contradicts to the second identity in (\ref{2.8}). This yields $\hat{i}(y_2)>0$.

Now by Lemma 4.6, when the case of (\ref{4.49}) holds, both $\hat{i}(y_2)$ and $\hat{i}(y_2)$ must be irrational numbers
by the definitions of $\th_i$s. When the case of (\ref{4.50}) holds, $\hat{i}(y_2)$ must be an irrational number by the
definition of $\th_3$. Then $\hat{i}(y_1)$ must be irrational too by the first identity of (\ref{4.8}). \hfill\hb

\medskip

{\bf Step B:} {\it The aim in this Step B is to prove the following theorem.}

\medskip

{\bf Theorem 4.8.} {\it Let $\Sigma$ be a non-degenerate $C^3$ compact star-shaped hypersurface
in $\R^{6}$ possessing precisely $2$ prime closed characteristics $(\tau_1,y_1)$ and $(\tau_2,y_2)$.
Then their Maslov-type indices satisfy $\{i(y_1,1), i(y_2,1)\} = \{0, -1\}$. Consequently Theorem
1.1 holds. }

\medskip

{\bf Proof.} In fact, by assuming that $^{\#}\T(\Sg)=2$, we will obtain some precise information
(\ref{4.56})-(\ref{4.57}) in Lemma 4.9, and (\ref{4.67})-(\ref{4.68}) in Lemma 4.10 about the indices
of $y_j,j=1,2$, which imply $\{i(y_1,1),i(y_2,1)\}=\{0,-1\}$ via Theorem 3.4. Because it is assumed in Theorem 1.1 that $\Sigma$ does not possess any prime closed characteristic with either Maslov-type index $0$, or Maslov-type index $-1$, so Theorem 1.1 holds. \hfill\hb

Next we focus on the proofs of the following Lemma 4.9 and Lemma 4.10.

\medskip

{\bf Lemma 4.9.} {\it Assume that $^{\#}\T(\Sg)=2$ and the two prime closed characteristics $(\tau_1,y_1)$ and $(\tau_2,y_2)$
satisfy (\ref{4.49}), then we must have
\bea
&& i(y_1)=-3, \qquad i(y_1^m) = 2\left[\frac{m\th_1}{2\pi}\right] - 3, \qquad \forall\;m\in\N, \lb{4.56}\\
&& i(y_2) = -4,  \qquad i(y_2^m) = - m + 2\left[\frac{m\th_2}{2\pi}\right] - 3, \qquad \forall\;m\in\N. \lb{4.57}
\eea}

{\bf Proof.} Note first that both $\hat{i}(y_1)>0$ and $\hat{i}(y_2)>0$ hold by Lemma 4.7.

Note that by (4.49) and Theorems 8.1.4, 8.1.7 and 8.2.1 of \cite{Lon3}, we have $i(y_1,1)\in 2\Z$ and then
$i(y_1)\in 2\Z+1$ by Theorem 3.4. By Theorem 3.1 and Theorem 3.4 for every $m\in\N$ we have
\bea
i(y_1,m) &=& mi(y_1,1) + 2\left[\frac{m\th_1}{2\pi}\right],  \lb{4.58}\\
i(y_1^m) &=& m(i(y_1)+3) + 2\left[\frac{m\th_1}{2\pi}\right] - 3,  \lb{4.59}\\
\hat{i}(y_1) &=& i(y_1) + 3 + \frac{\th_1}{\pi} > 0.  \lb{4.60}\eea
Consequently from $\th_1\in (0,2\pi)$ we have $i(y_1) > -5$. Because it is odd, we obtain $i(y_1)\ge -3$, i.e,
\be  i(y_1) \in 2\N-5.  \lb{4.61}\ee

By (4.49) and Theorems 8.1.4, 8.1.7 and 8.2.1 of \cite{Lon3}, we have $i(y_2,1)\in 2\Z+1$ and then
$i(y_1)\in 2\Z$ by Theorem 3.4. By Theorem 3.1 and Theorem 3.4 for every $m\in\N$ we have
\bea
i(y_2,m) &=& mi(y_2,1) + 2\left[\frac{m\th_2}{2\pi}\right],  \lb{4.62}\\
i(y_2^m) &=& m(i(y_2)+3) + 2\left[\frac{m\th_2}{2\pi}\right] - 3,  \lb{4.63}\\
\hat{i}(y_2) &=& i(y_2) + 3 + \frac{\th_2}{\pi} > 0.  \lb{4.64}\eea
Consequently from $\th_j\in (0,2\pi)$ we have $i(y_2) > -5$. Because it is even, we obtain $i(y_2)\ge -4$, i.e.,
\be  i(y_2) \in 2\N - 6.  \lb{4.65}\ee

Then by (\ref{2.10}), we have $\hat{\chi}(y_1) = -1$ and $\hat{\chi}(y_2) = \frac{1}{2}$.
Thus the first identity in (\ref{2.9}) becomes
\be  - \frac{1}{i(y_1) + 3 + \aa_1} + \frac{1}{2(i(y_2) + 3 + \aa_2)} = \frac{1}{2}.  \lb{4.66}\ee
where we let $\aa_j=\frac{\th_j}{\pi}$ for $j =1, 2$.

We continue the proof in two cases.

{\bf Case 1.} {\it $i(y_1)\ge -3$ and $i(y_2)\ge -2$.}

In this case, because $\aa_j>0$ with $j=1, 2$, (\ref{4.66}) yields
$$  \frac{1}{2} > \frac{1}{2(1 + \aa_2 + 2 + i(y_2))} = \frac{1}{2} + \frac{1}{3 + i(y_1) + \aa_1} > \frac{1}{2}. $$
This contradiction shows that Case 1 can not happen.

{\bf Case 2.} {\it $i(y_1)\ge -1$ and $i(y_2) = -4$.}

By Lemma 2.2, (\ref{2.12}) and $i(y_j^m)\ge i(y_j)\ge -4$ for any $m\in\N$ and $j=1, 2$, there holds
$M_{-3}-M_{-4}\ge 0$. However, from Lemma 2.2, (\ref{4.59}) and (\ref{4.63}), it turns out that $M_{-3}=0, M_{-4}\ge1$
in this case. This contradiction shows that Case 2 can not happen.

\medskip

Thus by Case 1 and Case 2, it yields $i(y_1)=-3$ and $i(y_2)=-4$. So Lemma 4.9 is proved. \hfill\hb

\medskip

{\bf Lemma 4.10.} {\it Assume that $^{\#}\T(\Sg)=2$ and the two prime closed characteristics $(\tau_1,y_1)$ and $(\tau_2,y_2)$
satisfy (\ref{4.50}), then we must have
\bea
&& i(y_1) = -4,\quad i(y_1^m) = - 2m + 2\sum_{j=1}^2\left[\frac{m\th_j}{2\pi}\right] - 2, \qquad \forall\;m\in\N, \lb{4.67}\\
&& i(y_2)= -3,\quad i(y_2^m) = 2\left[\frac{m\th_3}{2\pi}\right] - 3, \qquad \forall\;m\in\N. \lb{4.68}\eea}

{\bf Proof.} Note first that both $\hat{i}(y_1)>0$ and $\hat{i}(y_2)>0$ hold by Lemma 4.7.

Note that by (4.50) and Theorems 8.1.4 and 8.1.7 of \cite{Lon3}, we have $i(y_1,1)\in 2\Z+1$ and then
$i(y_1)\in 2\Z$ by Theorem 3.4. By Theorem 3.1 and Theorem 3.4 for every $m\in\N$ we have also
\bea
i(y_1,m) &=& m(i(y,1)-1) + 2\sum_{j=1}^2\left[\frac{m\th_j}{2\pi}\right] +1,  \lb{4.69}\\
i(y_1^m) &=& m(i(y_1)+2) + 2\sum_{j=1}^2\left[\frac{m\th_j}{2\pi}\right] - 2,  \lb{4.70a}\\
\hat{i}(y_1) &=& i(y_1) + 2 + \frac{\th_1+\th_2}{\pi}.  \lb{4.71}\eea

So by $\th_1, \th_2\in (0,2\pi)$ it yields $i(y_1)>-6$. Because $i(y_1)$ is even, we obtain $i(y_1) \ge -4$ and
\bea i(y_1^m)=i(y_1)=0\ (\mod 2),\quad \forall m\ge1.\lb{4.72}\eea

By (4.50) and Theorems 8.1.4, 8.1.7 and 8.2.1 of \cite{Lon3}, we have $i(y_2,1)\in 2\Z$ and then
$i(y_2)\in 2\Z+1$ by Theorem 3.4. By Theorem 3.1 and Theorem 3.4 for every $m\in\N$ we have also
\bea
i(y_2,m) &=& mi(y_2,1) + 2\left[\frac{m\th_3}{2\pi}\right],  \lb{4.73}\\
i(y_2^m) &=& m(i(y_2)+3) + 2\left[\frac{m\th_3}{2\pi}\right] - 3,  \lb{4.74a}\\
\hat{i}(y_2) &=& i(y_2) + 3 + \frac{\th_3}{\pi}.  \lb{4.75}\eea

Thus by $\th_3\in (0,2\pi)$ it yields $i(y_2) > -5$. Because $i(y_2)$ is odd, we obtain
$i(y_2) \ge -3$ and
\bea i(y_2^m)=i(y_2)=1\ (\mod 2),\quad \forall m\ge1.\lb{4.76}\eea

Then by (\ref{2.10}), we have $\hat{\chi}(y_1) = 1$ and $\hat{\chi}(y_2) = -1$. Thus the first identity
in (\ref{2.9}) becomes
\be \frac{1}{i(y_1) + 2 + \aa_1} - \frac{1}{i(y_2) + 3 + \aa_2} = \frac{1}{2}, \lb{4.77}\ee
where we let $\aa_1=\frac{\th_1+\th_2}{\pi}$ and $\aa_2=\frac{\th_3}{\pi}$.

We continue the proof in three cases according to the values of $i(y_1)\ge -4$ and $i(y_2)\ge -3$.

{\bf Case 1.} {\it $i(y_1)\ge 0$ and $i(y_2)\ge -3$.}

In this case, (\ref{4.77}) yields
$$  \frac{1}{2} > \frac{1}{2 + i(y_1) + \aa_1} = \frac{1}{2} + \frac{1}{3 + i(y_2) + \aa_2} > \frac{1}{2}. $$
This contradiction shows that Case 1 can not happen.

{\bf Case 2.} {\it $i(y_1)=-2$ and $i(y_2)\ge -3$.}

It follows from Theorem 3.2 that there exists infinitely many 3-tuples $(m_{1},m_{2},N)$ such that
\bea
i(y_j,2m_j+1) &=& 2 N+i(y_j,1), \quad j=1,2,                        \\
i(y_j,2m_j-1) &=&  2 N-i(y_j,1)-2, \quad j=1,2,  \\
i(y_1, 2m_1)&=& 2 N -3+2\Delta_1,    \\
i(y_2, 2m_2)&=& 2 N -2+2\Delta_2,     \eea
which, together with Theorem 3.4, $\Delta_1\le 2$ and $\Delta_2\le 1$, implies that
\bea
i(y_{j}^{2m_j+1}) &=& 2 N+i(y_{j}), \quad j=1,2,                       \lb{CIJT1}  \\
i(y_{j}^{2m_j-1}) &=&  2 N-i(y_{j})-8, \quad j=1,2,  \lb{4.83}\\
i(y_{1}^{2m_1})&=& 2 N -6+2\Delta_1\le 2N-2,    \lb{4.84}\\
i(y_{2}^{2m_2})&=& 2 N -5+2\Delta_2\le 2N-3.     \lb{CIJT2}\eea

Note that, in this case, it follows from (\ref{4.70a}) and (\ref{4.74a}) that ${i(y_{j}^{m})}$ is an increasing
sequence on $m\in\N$ for $j=1,2$, respectively. Then by $i(y_1)=-2$ and $i(y_2)\ge -3$, it further yields
\bea
&&-2=i(y_1)\le i(y_{1}^{m}) \le i(y_{1}^{2m_1-1})=2 N-6,\qquad\forall\ 1\le m\le 2m_1-1,         \lb{4.86a}  \\
&&-3\le i(y_2)\le i(y_{2}^{m}) \le i(y_{2}^{2m_2-1})\le 2 N-5,\qquad\forall\ 1\le m\le 2m_2-1,  \lb{4.87a}\\
&&i(y_{2}^{m}) \ge i(y_{2}^{2m_2+1})\ge 2N-3,\qquad\forall\ m\ge 2m_2+1,   \lb{4.88a}\\
&&i(y_{1}^{m}) \ge i(y_{1}^{2m_1+1})\ge 2N-2,\qquad\forall\ m\ge 2m_1+1.    \lb{4.89a}
\eea

Then by Lemma 2.2, (\ref{4.72}), (\ref{4.76}) and (\ref{4.84})-(\ref{4.89a}), we obtain
\bea 2m_1-1 &\le& M_{2N-6}+M_{2N-8}+\cdots+M_{-2},\lb{4.90}\\
2m_2-1 &\le& M_{2N-5}+M_{2N-7}+\cdots+M_{-3}\le 2m_2,\lb{4.91}\\
2m_1 &\ge& M_{2N-4}+M_{2N-6}+\cdots+M_{-2}.\lb{4.92}
\eea

Then, together with the Morse inequality (\ref{2.12}), there holds
\bea 2m_2-(2m_1-1)&\ge& M_{2N-5}-M_{2N-6}+M_{2N-7}-M_{2N-8}+\cdots+M_{-3}-M_{-2}\nn\\
                  &\ge& b_{2N-5}-b_{2N-6}+b_{2N-7}-b_{2N-8}+\cdots+b_{-3}-b_{-2}\nn\\
                  &=&-N+2,\\
2m_1-(2m_2-1)&\ge& M_{2N-4}-M_{2N-5}+M_{2N-6}-M_{2N-7}+\cdots+M_{-2}\nn\\
                  &\ge& b_{2N-4}-b_{2N-5}+b_{2N-6}-b_{2N-7}+\cdots+b_{-2}\nn\\
                  &=&N-1,\lb{4.94}
\eea
which implies $-2\le 2m_1-2m_2-N\le -1$. However, by Theorem 3.3, we can choose 3-tuple
$(m_{1},m_{2},N)$ such that $(4m_{1},4m_{2},4N)$ also satisfy (\ref{CIJT1})-(\ref{CIJT2}). And
by the same argument as those in (\ref{4.90})-(\ref{4.94}), there holds $-2\le 8m_1-8m_2-4N\le -1$,
which is impossible.

This contradiction shows that Case 2 can not happen.

{\bf Case 3.} {\it $i(y_1)=-4$ and $i(y_2)\ge -1$.}

By Lemma 2.2, (\ref{2.12}) and $i(y_j^m)\ge i(y_j)\ge -4$ for any $m\in\N$ and $j=1, 2$, there
holds $M_{-3}-M_{-4}\ge 0$. However, from Lemma 2.2, (\ref{4.70a}) and (\ref{4.74a}), it turns out
that $M_{-3}=0, M_{-4}\ge1$ in this case. This contradiction shows that Case 3 can not happen.

\medskip

Thus by Cases 1-3 it yields $i(y_1)=-4$ and $i(y_2)=-3$. The proof of Lemma 4.10 is finished. \hfill\hb

\bibliographystyle{abbrv}

\end{document}